\documentclass{article}

\usepackage{graphicx}
\usepackage{amssymb}
\usepackage{amsthm}
\usepackage{pstricks,pst-node,pst-coil,pst-plot}
\usepackage{enumerate}

\newdimen\templaenge

\newtheorem{theorem}{Theorem}[section]
\newtheorem{lemma}[theorem]{Lemma}
\newtheorem{proposition}[theorem]{Proposition}
\newtheorem{corollary}[theorem]{Corollary}
\theoremstyle{definition}
\newtheorem*{definition}{Definition}
\newtheorem*{remark}{Remark}
\newtheorem*{example}{Example}
\let\al\alpha
\let\be\beta
\let\ga\gamma
\let\de\delta
\let\ep\varepsilon

\let\la\lambda
\let\om\omega
\let\ph\varphi
\let\si\sigma
\let\th\theta
\let\Ga\Gamma
\let\Si\Sigma
\let\Th\Theta

\let\na\nabla
\let\pa\partial

\def\rear{\mbox{$\Rightarrow$}}

\def\area{\mathrm{area}}

\def\ghyp{g_{\mathrm{hyp}}}
\def\qarf{q_\mathrm{spin}}
\def\Arf{\mathrm{Arf}}
\def\Hom{\mathrm{Hom}}
\def\ie{i.\thinspace e.\ }
\def\eref#1{(\ref{#1})}
\let\lan\langle
\let\ran\rangle

\def\res#1#2{{#1}\lower .11ex\hbox{$|$}\lower .644ex\hbox{$\scriptstyle #2$}}
\def\proof#1{{\par\medbreak\noindent {\bf Proof\setbox0\hbox{#1}%
\ifdim\wd0=0pt .\else\ \ignorespaces #1.\fi}\enspace}}
\def\iop#1{{\par\medbreak\noindent {\bf Idea of proof\setbox0\hbox{#1}%
\ifdim\wd0=0pt .\else\ \ignorespaces #1.\fi}\enspace}}

\def\fm{T_{-k}}
\def\fp{T_{k}}
\def\fmp{\fm\cup \fp}
\def\stabnorm#1{{\|#1\|_{\rm st}}}
\def\length#1{\mathrm{length}(#1)}
\def\SO{{\mathop{\rm SO}}}
\def\Spin{{\mathop{\rm Spin}}}

\let\ti\tilde
\let\witi\widetilde
\def\ccutM{\widetilde{M}}

\def\spinsys{{\mathop{\textrm{spin-sys}}}}

\newcommand{\nc}{\newcommand}
\nc{\CC}{\mathbb C}
\nc{\HH}{\mathbb H}
\nc{\NN}{\mathbb N}
\nc{\RR}{\mathbb R}
\nc{\ZZ}{\mathbb Z}

\def\qedmath#1{\setbox0\hbox{$\displaystyle #1$}\templaenge=\textwidth\advance\templaenge by -\wd0%
\setbox1\hbox{$\Box$}\advance\templaenge by -2\wd1%
$$#1\hbox to0pt{\kern.5\templaenge$\Box$\kern-.5pt\hss}$$\par\bigbreak}         

\long\def\komment#1{}

\def\lss{\mbox{spec}_{L^2}}
\def\barM{\overline{M}}
\def\vol{\mbox{vol}}

\begin{document}

\title{Dirac eigenvalue estimates on surfaces}

\author{Bernd Ammann and Christian B\"ar
\footnote{Both authors partially supported 
by the European Contract Human Potential Programme,
Research Training Networks HPRN-CT-2000-00101 and HPRN-CT-1999-00118}
}

\date{14. May 2001}
\maketitle

\begin{abstract}
\noindent
We prove lower Dirac eigenvalue bounds for closed surfaces with a spin 
structure whose Arf invariant equals 1. 
Besides the area only one geometric quantity enters in these estimates, 
the spin-cut-diameter $\de(M)$ which depends on the choice of spin structure.
It can be expressed in terms of various distances on the surfaces or,
alternatively, by stable norms of certain cohomology classes.
In case of the 2-torus we obtain a positive lower bound for all Riemannian
metrics and all nontrivial spin structures.
For higher genus $g$ the estimate is given by
$$
|\la| \geq {2\sqrt{\pi}\over (2g+1)\,\sqrt{\area(M)}} -{1\over \de(M)}.
$$
The corresponding estimate also holds for the $L^2$-spectrum of the Dirac
operator on a noncompact complete surface of finite area.
As a corollary we get positive lower bounds on the Willmore integral
for all 2-tori embedded in $\RR^3$.\\


\end{abstract}

%%%%%%%%%%%%%%%%%%%%%%%%%%%%%%%%%%%
\section{Introduction}
%%%%%%%%%%%%%%%%%%%%%%%%%%%%%%%%%%%

Relating analytic invariants of the Dirac operator such as the
eigenvalues to the geometry of the underlying manifold is in general
a difficult problem.
Explicit computation of the spectrum is possible only in cases
of very large symmetry, see 
\cite{baer:92a},\cite{baer:94},\cite{baer:96},\cite{bunke:91},\cite{cahen.franc.gutt:89,cahen.franc.gutt:94},\cite{camporesi.higuchi:96},\cite{fegan:87},\cite{friedrich:84},\cite{hitchin:74},\cite{milhorat:p},\cite{milhorat:92},\cite{pfaeffle:p99},\cite{seegerdiplom},\cite{seeger:99},\cite{seifarth.semmelmann:93},\cite{strese:80b},\cite{strese:80a},\cite{sulanke:79},\cite{trautman:93,trautman:95} for examples.
In general, the best one can hope for are geometric bounds on the eigenvalues.
The first lower eigenvalue bounds
\cite{friedrich:80}, \cite{kirchberg:86},\cite{kirchberg:88},\cite{kramer.semmelmann.weingart:98},\cite{kramer.semmelmann.weingart:99} for the Dirac spectrum require
positivity of the scalar curvature since they are based on variations of the 
Lichnerowicz formula $D^2=\nabla^\ast\nabla + \mbox{scal}/4$.
Refining this technique Hijazi \cite{hijazi:86,hijazi:91}
could estimate the smallest Dirac eigenvalue against the corresponding
eigenvalue of the Yamabe operator.
A completely different approach building on Sobolev embedding theorems
was used by Lott \cite{lott:86} and the first author \cite{ammann:p00b} 
to show that
for each closed spin manifold $M$ and each conformal class $[g_0]$ on $M$
there exists a constant $C=C(M,[g_0])$ such that all nonzero Dirac 
eigenvalues $\la$ with respect to all Riemannian metrics $g\in [g_0]$ satisfy
$$
\la^2\geq {C\over\vol(M)^{2/n}}.
$$ 
On the 2-sphere $M=S^2$ there is only one conformal class of metrics
(up to the action of the diffeomorphism group) and we therefore get
a nontrivial lower bound for all metrics.
Lott conjectured that in this case the optimal constant should be 
$C=4\pi$.
Returning to the Bochner technique the second author showed that this 
is in fact true:

\sloppy

\begin{theorem}[{\cite[Theorem 2]{baer:92b}}]\label{baersphere}
Let $\la$ be any Dirac eigenvalue of the $2$-sphere $S^2$ equipped with an 
arbitrary Riemannian metric. 
Then  
$$
\la^2\geq {4\pi\over\area(S^2)}.
$$ 
Equality is attained if and only if $S^2$ carries a metric of 
constant Gauss curvature.
\end{theorem}

In particular, there are no harmonic spinors on $S^2$.
Theorem~\ref{baersphere} will be the central tool to derive our new estimates
in the present paper.
Examples \cite{baer:96}, \cite{seegerdiss} show that such an estimate 
is neither possible for higher dimensional spheres nor for surfaces of 
higher genus, at least not in this generality.
Every closed surface of genus at least 1 has a spin structure and a metric
such that 0 is an eigenvalue, i.~e.\ there are nontrivial harmonic spinors
\cite{friedrich:84}, \cite{hitchin:74}.
The 2-torus $T^2$ has four spin structures one of which is called trivial and
the others nontrivial.
Provided with the trivial spin structure, $T^2$ has harmonic spinors for all
Riemannian metrics.
On the other hand, for the three nontrivial spin structures 0 is never
an eigenvalue.
So it should in principle be possible to give a geometric lower bound
in this latter case.
The problem is that this estimate must take into account the choice of 
spin structure but the Bochner technique is based on local computation
where the spin structure is invisible.
Hence new techniques are needed.

The first estimate using information from the choice of spin structure
has been derived by the first author \cite[Corollary~2.4]{ammann:p00a}. 
On a torus with a Riemannian metric and a nontrivial spin structure 
there is a lower bound for any eigenvalue $\la$ of the Dirac operator.
Let $K$ denote Gauss curvature.
Recall that the systole is the minimum of the lengths of all noncontractible 
closed curves. 
The spinning systole $\spinsys(T^2)$ is the minimum of the lengths of all 
noncontractible simple closed curves, along which the spin structure is 
nontrivial. If there exists $p>1$ with $\|K\|_{L^p}\cdot
\area(T^2)^{1-(1/p)}<4\pi$, then there is a positive number $C>0$ such that
  $$
  \la^2 \geq {C\over \spinsys(T^2)^2}.
  $$
Here $C$ is an explicitly given expression in $p$, $\|K\|_{L^p}$, 
the area, and the systole.

The Arf invariant associates to each spin structure on a closed surface
the number $1$ or $-1$.
In case of the 2-torus the Arf invariant of the trivial spin structure
is $-1$ while the three nontrivial spin structures have Arf invariant $1$.
In the present paper we prove explicit geometric lower bounds for the first 
eigenvalue of the square of the Dirac operator on closed 
surfaces $M$ of genus $\geq 1$ provided the spin structure has Arf 
invariant $1$.
Only two geometric quantities enter, the area of the surface and an
invariant we call the spin-cut-diameter $\de(M)$.
The number $\de(M)$ is defined by looking at distances between loops in the 
surface along which the spin structure is nontrivial and which are 
linearly independent in homology. 
It exists if and only if the Arf invariant of the spin structure equals~$1$.
It can also be defined in terms of stable norms of certain cohomology 
classes which depend on the choice of spin structure 
(Proposition~\ref{dusys}).

In the case of a 2-torus we show:

{\bf Theorem~\ref{2torus}.}
{\em Let $T^2$ be the 2-torus equipped with an arbitrary Riemannian metric and 
a spin structure whose Arf invariant equals $1$. 
Let $\la$ be an eigenvalue of the Dirac operator and let $\de(T^2)$ be
the spin-cut-diameter. 
Then for any $k\in \NN$, 
$$
|\la|\geq -{2\over k\,\de(T^2)} + \sqrt{{\pi\over k\,\area(T^2)}+{2\over k^2 \de(T^2)^2}}.
$$
}

The right hand side of this inequality is positive for sufficiently 
large $k$.
Hence this theorem gives a nontrivial lower eigenvalue bound for the 
Dirac operator for all Riemannian metrics and all nontrivial spin structures
on the 2-torus.

Similarly, for higher genus we obtain:

{\bf Theorem~\ref{ghoch}.}
{\em Let $M$ be a closed surface of genus $g\ge 1$ with a Riemannian metric 
and a spin structure whose Arf invariant equals 1.
Let $\de(M)$ be the spin-cut-diameter of $M$.
Then for all eigenvalues $\la$ of the Dirac operator we have
  $$|\la| \geq {2\sqrt{\pi}\over (2g+1)\,\sqrt{\area(M)}} -{1\over \de(M)}.$$
}

In the case $g=1$ this estimate is simpler but weaker than Theorem~\ref{2torus}.
Every surface of genus $g\geq 2$ admits metrics and spin structures
such that this estimate is nontrivial. 
But in contrast to the first theorem there are also Riemannian
metrics and spin structures on surfaces of genus $g\geq 1$ for which the 
right hand side of this inequality is negative although there are no 
harmonic spinors.

If one restricts one's attention to surfaces embedded in $\RR^3$, then
one has the Willmore integral $W(M)$ defined as the integral of the square of
the mean curvature.
It is well-known that the Willmore integral can be estimated against
Dirac eigenvalues.
Thus as a corollary to Theorem~\ref{2torus} we obtain

{\bf Theorem~\ref{willmore}}
{\em Let $T^2\subset\RR^3$ be an embedded torus.
Let $\de(T^2)$ be its spin-cut-diameter and let $W(T^2)$ be its Willmore
integral.
Then for any $k\in \NN$
$$
\sqrt{W(T^2)} \ge \sqrt{\frac{\pi}{k}+\frac{2\,\area(T^2)}{k^2\,\de(T^2)^2}}
- \frac{2\sqrt{\area(T^2)}}{k\,\de(T^2)}
$$
}

In the end of the paper we show that our spectral estimates also work
for noncompact complete surfaces of finite area.
In this case the spectrum need not consist of eigenvalues only.
We estimate the fundamental tone of the square of the Dirac operator
which gives the length of the spectral gap about $0$ in the $L^2$-spectrum,
see Theorem~\ref{ghochnc}.

The paper is organized as follows.
We start by recalling some basic definitions related to spin structures
and Dirac operators on surfaces.
We put some emphasis on the case of a surface embedded in $\RR^3$.
We then recall the Arf invariant and define the  spin-cut-diameter
$\de(M)$.
In Section~\ref{stable} we show how $\de(M)$ relates to the stable norm
of certain cohomology classes.
In Sections~\ref{torussec} and \ref{hochsec} we prove Theorems~\ref{2torus} 
and \ref{ghoch}.
The central idea of proof consists of constructing a surface of genus $0$
out of the given surface by cutting and pasting.
Then we apply Theorem~\ref{baersphere}.
The estimate for the Willmore integral is proved in 
Section~\ref{willmoresec} and in Section~\ref{noncompact} we study the 
$L^2$-spectrum of noncompact complete surfaces of finite area.

%%%%%%%%%%%%%%%%%%%%%%%%%%%%%%%%%%%
\section{Dirac operators on surfaces}
\label{diracsurfsec}
%%%%%%%%%%%%%%%%%%%%%%%%%%%%%%%%%%%

Let $M$ be an oriented surface with a Riemannian metric. 
Rotation by $90$ degrees in the positive direction defines a complex 
multiplication $J$ on $TM$.
The bundle $\SO(M)$ of oriented orthonormal frames is an 
$S^1$-principal bundle over $M$.
Let $SM$ be the bundle of unit tangent vectors on $M$.
Then  $v\mapsto (v,Jv)$ is a fiber preserving diffeomorphism from $SM$ 
to $\SO(M)$ with inverse given by projection to the first vector.

Let $\Th:S^1\to S^1$ be the nontrivial double covering of $S^1$. 
A spin structure on $M$ is an $S^1$-principal bundle $\Spin(M)$ over $M$
together with a twofold covering map $\th:\Spin(M)\to \SO(M)$ such that the
diagram
\begin{equation}\label{spincompat}
\begin{array}{cccl}
\Spin(M) \times S^1 &\rightarrow& \Spin(M) & \\
& & &\searrow \\
\downarrow\th\times\Theta & &   \downarrow\th & \quad M\\
& & &\nearrow \\
\SO(M) \times S^1 &\rightarrow& \SO(M) &
\end{array}
\end{equation}
commutes.

Every orientable surface admits a spin structure, but it is in general
not unique. 
The number of possible
spin structures on $M$ equals the number of elements in $H^1(M,\ZZ_2)$.

\begin{example}
Let $i:M\hookrightarrow \RR^3$ be an immersion of an oriented surface (not necessarily compact, and
possibly with boundary) into $\RR^3$. We define a map $i_*:\SO(M)\to \SO(3)$ as follows:
$(v,Jv)\in\SO(M)$ over a basepoint $m\in M$ is mapped to 
$(v,Jv,v\times Jv)\in \SO(3)$.
Here $\times$ denotes the vector cross product in $\RR^3$.
Let $\Spin(M)$ be the pullback of the double covering $\Theta_3:\Spin(3)\to \SO(3)$,
\ie  
$$
\Spin(M):=\left\{\big((v,Jv),A\big)\in\SO(M)\times \Spin(3)\,\Big|\,i_*(\SO(M))=\Theta_3(A)\right\}.
$$
Then $\Spin(M)\to\SO(M)$ is a fiberwise nontrivial double covering.
Let $\pi:\SO(M)\times \Spin(3)\to\SO(M)$ be the projection onto 
the first component. Then $(\Spin(M),\res{\pi}{\Spin(M)})$ is a spin structure on $M$,
the \emph{spin structure induced by the immersion}. 
\end{example}

Let $\ga:S^1\to M$ be an immersion or, in other words, a regular closed curve.
Then the vector field $\dot\ga\over |\dot\ga|$ is a section of $SM$ 
along $\ga$, which, by the above diffeomorphism from $SM$ to $\SO(M)$,
yields the section $({\dot\ga\over |\dot\ga|},J{\dot\ga\over |\dot\ga|})$
of $\SO(M)$ along $\ga$. 

\begin{definition}
The spin structure $(\Spin(M),\th)$ is said to be \emph{trivial along} $\ga$ 
if this section lifts to a closed curve in $\Spin(M)$ via $\th$. 
\end{definition}

This notion is invariant under homotopic deformation of $\ga$ within the 
class of immersions.

\begin{example}
The unique spin structure on $\RR^2$ is nontrivial along any simple closed curve.
More generally, any spin structure on a surface $M$ is nontrivial along any 
contractible simple closed curve.
\end{example}

\begin{proposition}\label{nontrivcrit}
Let $i:M\hookrightarrow\RR^3$ be an immersion. Let $\ga:S^1\to M$ be a simple closed curve.
If $\ga$ is a parametrization of the boundary of an immersed two-dimensional 
disk $j:D\hookrightarrow\RR^3$ intersecting $i(M)$ transversally, then the 
spin structure on $M$ induced by $i$ is nontrivial
along $\ga$.  
\end{proposition}

%%%%%
%%
%% Bild: Disk in Torus
%% former: diskintorus.pic
%%
%%%%%

\begin{figure}[h]
\begin{center}\includegraphics*[width=10cm]{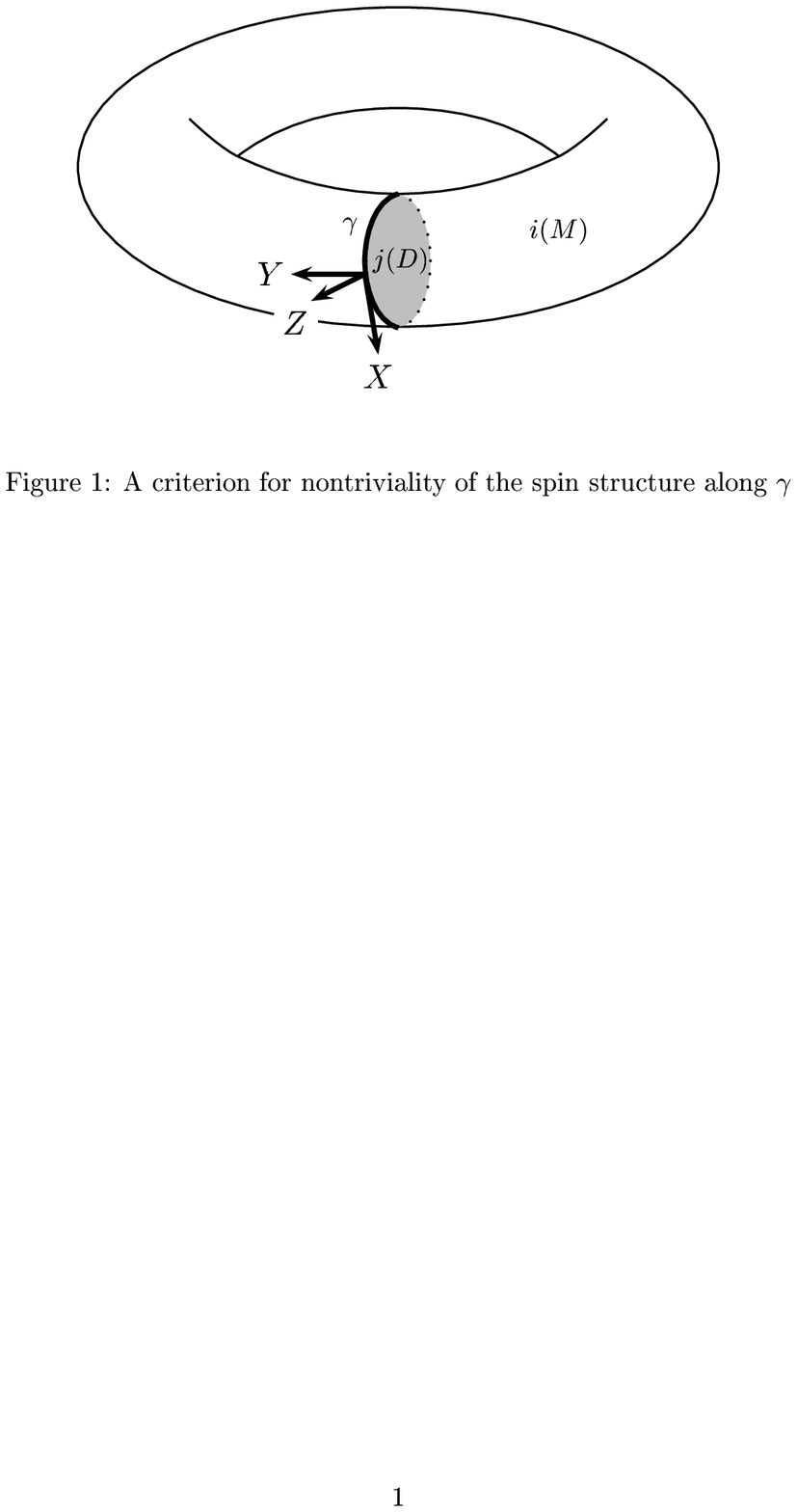}\end{center}
\end{figure}
\stepcounter{figure}

\proof{}
We can assume that $j(D)$ and $i(M)$ intersect orthogonally along $\ga$.
We set $X(t):={\dot\ga(t)\over |\dot\ga(t)|}$, $Y(t):=J_M X(t)$ and 
$Z(t):=X(t)\times Y(t)$.
The induced spin structure on $M$ is trivial along $\ga$ if and only if 
$S^1\to \SO(3), t\mapsto (X(t),Y(t),Z(t))$, lifts 
to a closed loop in $\Spin(3)$.
Analogously, we view $\ga$ as a curve on $j(D)$, 
we define the vector fields $\hat Y(t):=J_D X(t)$ 
and $\hat Z(t):=X(t)\times \hat Y (t)$. Because of the orthogonality of 
$j(D)$ and $i(M)$ we have $\hat Y (t)=\pm Z(t)$ and  $\hat Z(t) = \mp Y(t)$.

Hence $t\mapsto (X(t),\hat Y(t),\hat Z(t))$ lifts
to $\Spin(3)$ if and only if $t\mapsto (X(t),Y(t),Z(t))$ lifts.
The induced spin structure on $M$ is nontrivial along $\ga$ if and only 
if the spin structure on $D$ is nontrivial along $\ga$.
This is always true according to the previous example.
\qed

\begin{example}
Let $Z:=\{(x,y,z)\in\RR^3\,|\,x^2+y^2=1\}$ be the cylinder with the 
induced spin structure. Let $\ga$ be any simple closed curve in $Z$.
We show that the spin structure is nontrivial along $\ga$: 
If $\ga$ is contractible then the spin structure is nontrivial because of
the preceeding example. 
If $\ga$ is noncontractible, then
$[\ga]$ generates $\pi_1(Z)$ (Lemma~\ref{appendixlemmaA}). 
Hence it bounds a disk transversal to~$Z$.
\end{example}

Let $\Si^+M:=\Spin(M)\times_{\iota}\CC$ be the complex 
line bundle over $M$ associated to the $S^1$-principal bundle 
$\Spin(M)$ and to the standard representation $\iota:S^1\to U(1)$. 
This line bundle is called the \emph{bundle of positive half-spinors}, its complex 
conjugate $\Si^-M:=\overline{\Si^+M}$ is the 
\emph{bundle of negative half-spinors} and their sum 
$\Si M:=\Si^+M \oplus \Si^- M$ is the \emph{spinor bundle}.

Clifford multiplication consists of complex linear maps
\begin{eqnarray*}
  TM\otimes_\CC \Si^+ M & \to & \Si^- M\\
  \overline{TM}\otimes_\CC \Si^- M & \to & \Si^+ M
\end{eqnarray*}
denoted by $v\otimes \si\mapsto v\cdot \si$. It satisfies the 
Clifford relations
  $$v\cdot w \cdot \si + w\cdot v \cdot \si + 2\lan v,w\ran\, \si=0$$
for all $v,w\in TM$ and $\si \in \Si M$ over a common base point. 

The Levi-Civita connection on $TM$ gives rise to a connection-1-form on 
$\Spin(M)$ and this in turn defines a Hermitian connection $\na$ on $\Si M$.

\begin{definition}
The \emph{Dirac operator} $D$ is a map from smooth sections of $\Si M$ 
to smooth sections of $\Si M$ which is locally given by the formula
  $$D\Psi:=e_1\cdot\na_{e_1}\Psi + e_2\cdot\na_{e_2}\Psi$$
for a local orthonormal frame $(e_1,e_2)$ of $TM$.
\end{definition}

It is easily checked that the definition does not depend on the choice of 
the local frame and that $D$ is a formally self-adjoint elliptic
operator. Hence, if $M$ is closed, the spectrum of $D$ is real and 
discrete with finite multiplicities. 

For any smooth function $f$ and smooth spinor $\Psi$ the equation
$$
D(f\Psi) = \nabla f \cdot \psi + f D\Psi
$$
holds.
Here $\nabla f$ denotes the gradient of $f$.

For more background material on Dirac operators and spin structures see 
e.~g.\ \cite{lawson.michelsohn:89}, \cite{friedrich:buch}, or \cite{roe:88}.

To simplify notation a {\em closed} surface will always mean a surface
which is compact, without boundary, and {\em connected}.

%%%%%%%%%%%%%%%%%%%%%%%%%%%%%%%%%%%%%%%%%%%%%
\section{Arf invariant and spin-cuts}
\label{arfinvsection}\label{spincutsection}
%%%%%%%%%%%%%%%%%%%%%%%%%%%%%%%%%%%%%%%%%%%%%

In this section we review some properties of the Arf invariant 
which is an invariant of a spin structure 
on a surface (see \cite{pinkall:85a} for more details).
For closed oriented surfaces with spin structures whose Arf invariant 
equals~$1$ we define
a geometric quantity, the  spin-cut-diameter,
which will play an important role in our estimate.

Let $V$ be a $2g$-dimensional vector space over the field $\ZZ_2$, $g\in \NN$, together with a symplectic
$2$-form $\om:V\to \ZZ_2=\{0,1\}$. A \emph{quadratic form} on $(V,\om)$ is a map
$q:V\to \ZZ_2$, such that
  $$q(a+b) = q(a)+q(b)+\om(a,b)\quad \qquad a,b\in V.$$ 

The difference of two quadratic forms on $(V,\om)$ is a linear
map from $V$ to $\ZZ_2$ and vice versa the sum of a linear map 
$V\to\ZZ_2$ and a quadratic form is again a quadratic form.
Hence the space of quadratic forms on $V$ is an affine space over $\Hom(V,\ZZ_2)$.

\begin{example}
Let $M$ be a closed oriented surface.
Let $V:=H_1(M,\ZZ_2)$ and let $\om$ be the intersection form $\cap$.
Fix a spin structure on $M$. 
We associate to each spin structure a quadratic form $\qarf$ on $(V,\om)$ 
as follows. 
Each homology class $a\in H_1(M,\ZZ_2)$ is represented by an embedding 
$\ga:S^1\to M$.
We set $\qarf(a):=1$, if 
$(\dot \ga,J(\dot \ga)):S^1\to \SO(M)$ lifts to $\Spin(M)$, otherwise we 
set $\qarf(a):=0$.
\end{example}

According to Theorem~1 of \cite{kusner.schmitt:p96} the map 
$\qarf$ is a well-defined quadratic form on $(H_1(M,\ZZ_2),\cap)$.

The set of all spin structures on $M$ is an affine space over 
$H^1(M,\ZZ_2)=\Hom(H_1(M,\ZZ_2),\ZZ_2)$ and it is a well known
fact that the map which associates to any spin structure the corresponding 
quadratic form $\qarf$ is an isomorphism of affine $H^1(M,\ZZ_2)$-spaces 
from the space of spin structures on $M$ to the space of quadratic forms
on $(V,\om)=(H_1(M,\ZZ_2),\cap)$. 

\begin{definition}
For any quadratic form $q$ on $(V,\om)$ the \emph{Arf invariant} is defined by
\label{arfdef}
  $$\Arf(q) :={1\over {\sqrt{\# V}}}\sum_{a\in V} (-1)^{q(a)}.$$
The Arf invariant of a quadratic form corresponding to a spin structure 
will be called the \emph{Arf invariant} of that spin structure.
\end{definition}

\begin{lemma}\label{qsum}
Let $q_i$ be a quadratic form on $(V_i,\om_i)$ for $i=1,2$.
Then $q_1\oplus q_2$, given by
  $$(q_1\oplus q_2)(v_1+ v_2)=q_1(v_1) + q(v_2),$$
is a quadratic form on $(V_1\oplus V_2, \om_1\oplus\om_2)$. 
Moreover, 
  $$\Arf(q_1\oplus q_2)=\Arf(q_1)\Arf(q_2).$$
\end{lemma}
The proof is a simple counting argument.\qed

Any $2g$-dimensional symplectic vector space $V$ with a symplectic form $\om$ 
is isomorphic to the $g$-fold sum $V_2\oplus \cdots\oplus V_2$
where $V_2$ is the standard $2$-dimensional symplectic vector space.
Since the Arf invariants of the four possible choices of 
quadratic forms on $V_2$ are either $1$ or $-1$ the above lemma implies
  $$\Arf(q)\in\{-1,+1\}$$
for any quadratic form $q$ on any symplectic $\ZZ_2$-vector space.

\begin{proposition}\label{qaltern}
Let $q$ be a quadratic form on $(V,\om)$, $\dim V=2g$. Then the following statements are equivalent:
\begin{enumerate}[{\rm (1)}]
\item $\Arf(q)=1$.
\item There is a basis $e_1,f_1,\dots,e_g,f_g$ of\/ $V$ such that 
$\om(e_i,e_j)=\om(f_i,f_j)=0$, $\om(e_i,f_j)=\de_{ij}$, and 
$q(e_i)=q(f_j)=0$ for all $i,j$.
\item There are linearly independent vectors $e_1,\dots,e_g$ in $V$ such that 
${\om(e_i,e_j)=0}$ and $q(e_i)=0$ for all $i,j$.
\end{enumerate}
\end{proposition}

\proof{}
(2)\rear (1) follows directly from Lemma~\ref{qsum}. 

To show (3)\rear (2) let $e_1,\dots,e_g$ be linearly independent
vectors with $\om(e_i,e_j)=0$ and $q(e_i)=0$ for all $i,j$. 
Since $\om$ is symplectic, we can find $\tilde f_1,\dots,\tilde f_g$
satisfying $\om(e_i,\tilde f_j)=\de_{ij}$ and  
$\om(\tilde f_i,\tilde f_j)=0$ for all $i,j$.
If $q(\tilde f_i)=0$, we set $f_i:=\tilde f_i$, otherwise we put 
$f_i:=\tilde f_i+e_i$.

To see (1)\rear (3), we take a basis $e_1,f_1,\dots,e_g,f_g$ of $V$ satisfying
${\om(e_i,f_j)=\de_{ij}}$ and $\om(e_i,e_j)=\om(f_i,f_j)=0$. 
For every $i$ exactly one of the following holds: 
\begin{enumerate}[(a)]
\item $q(e_i)=q(f_i)=q(e_i+f_i)=1,\qquad\qquad$ or 
\item $q$ takes the value $0$ at exactly two of the vectors $e_i$, $f_i$ and 
$e_i+f_i$.
\end{enumerate}
In the second case, we can assume without loss of generality that $q(e_i)=q(f_i)=0$.
Let $I$ be the set of all $i$ for which (a) holds. Then by Lemma~\ref{qsum} $\Arf(q)=(-1)^{\# I}$.
If (1) holds, then $\# I$ is even, hence we may assume $I:=\{1,\dots,2k\}$.
For $j=1,\dots,k$ we replace $e_{2j-1}$ by $e_{2j-1}+f_{2j}$ and  $e_{2j}$ by $e_{2j}+f_{2j-1}$.
Then (3) holds.
\qed

\begin{example}
Let $M\hookrightarrow \RR^3$ be an embedded closed surface with the induced spin structure.
Then because of Propositions~\ref{nontrivcrit} and \ref{qaltern} (3) 
the Arf invariant of the spin structure is $1$.
As a consequence any immersion $M\hookrightarrow \RR^3$ whose induced spin structure has 
Arf invariant $-1$ is not regularly homotopic to an embedding.
\end{example}

\begin{remark}
In the literatur the $3$ spin structures on the 2-torus $T^2$ 
with Arf invariant  $1$
are called \emph{nontrivial} spin structures and the unique 
spin structure with Arf invariant $-1$ is called the \emph{trivial} 
spin structure.
\end{remark}

\begin{definition}
Let $M$ be a closed oriented surface of genus $g$. 
A \emph{cut} of $M$ is a family of pairwise disjoint simple closed 
curves $\ga_i:S^1\to M$, $i=1,\dots, g$, 
such that $[\ga_1],\dots,[\ga_g]$ are linearly independent in $H_1(M,\ZZ)$.
If, in addition, $M$ carries a spin structure, and if the spin structure is 
nontrivial along each of
the $\ga_i$, then we call $\ga_1,\ldots,\ga_g$ a \emph{spin-cut} of $M$.
\end{definition}

\begin{corollary}
Let $M$ be a closed oriented surface equipped with a spin structure. 
Then $M$ admits a spin-cut
if and only if the Arf invariant of the spin structure equals $1$.
\end{corollary}

\proof{}
If the Arf invariant is $1$, we can find vectors 
$e_1,\dots,e_g\in H_1(M,\ZZ_2)$ for which (3) of 
Proposition~\ref{qaltern} holds. 
For each $e_i$ we choose a preimage $\tilde e_i\in H_1(M,\ZZ)$
under the natural map $H_1(M,\ZZ)\to H_1(M,\ZZ_2)$. 
We choose $\tilde e_i$ such that $\tilde e_i$ is primitive,
\ie there are no $a_i\in H_1(M,\ZZ)$, $n\geq 2$ with $e_i=n \cdot a_i$. 
This choice can be made such that 
$\tilde e_i\cap \tilde e_j=0$ for all $i,j$. 
We choose a hyperbolic metric $\ghyp$ on $M$ and
represent $\tilde e_i$ by closed curves $\ga_i$ of minimal length. 
Then the $\ga_i$ are closed geodesics. 
They are simple closed curves because the $\tilde e_i$ are primitive.
Since $\tilde e_i\cap \tilde e_j=0$ and $\ghyp$ is hyperbolic,
$\ga_i$ and $\ga_j$ are disjoint for $i\neq j$. 
The spin structure is nontrivial along each $\ga_i$ because of $\qarf(e_i)=0$. 
Hence $\ga_1,\dots,\ga_g$ form a spin-cut of $M$.

Conversely, if $\ga_1,\dots,\ga_g$ form a spin-cut of $M$, then 
$[\ga_1],\dots,[\ga_g]\in H_1(M,\ZZ)$ form a linearly independent set 
of primitive elements in $H_1(M,\ZZ)$. 
Hence their images $e_i$ in $H_1(M,\ZZ_2)$ are also linearly independent. 
The $e_i$ satisfy (3) of Proposition~\ref{qaltern} and thus the Arf 
invariant is $1$.
\qed

\begin{definition}
Let $M$ be a closed surface.
Let $\ga_1,\dots,\ga_g$ be a cut. The \emph{cut-open} $\witi M$ of $M$ is a surface with boundary, such that there
is a smooth map $\witi M\to M$ which is a diffeomorphism from the interior of $\witi M$ onto $M\setminus \bigcup_{j=1}^g \ga_j$ and a 
twofold covering from the boundary $\pa \witi M$ 
onto $\bigcup_{j=1}^g \ga_j$. 
\end{definition}

%%%%%
%%
%% picture cutopen
%%
%%%%%

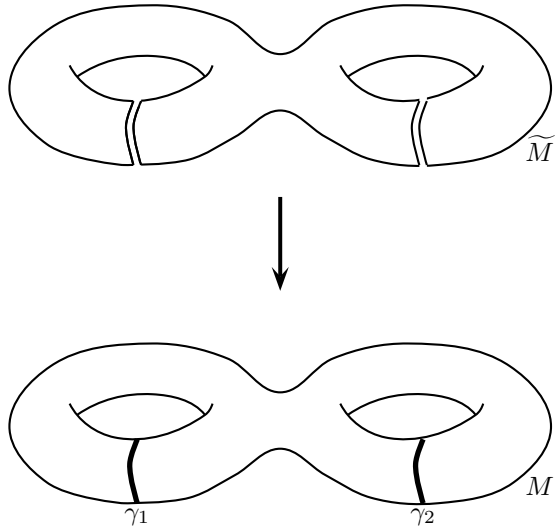
\begin{figure}[h]
\begin{center}
\psset{unit=5mm}
\begin{pspicture}(1,1)(14,15)
%   \psset{showpoints=true}
%    \psgrid(-2,-2)(13,20)

\def\Linse{
    \psecurve(1.8,5)(1.8,4.0)(3.0,3.1)(5.4,3.7)(5.6,4.0)(5.6,5)
    \psecurve(1.5,3)(2.0,3.7)(4,4.25)(5.4,3.7)(5.5,3)
    }

%% M unten
    \pscurve(3.2,1.35)(5,1.5)(6.2,2.1)(7.4,2.8)(8.6,2.1)(9.8,1.5)(11.6,1.35)%
    (13.6,1.9)(14.6,3.4)(13.5,5.0)(10.8,5.6)(8.8,5.2)(7.4,4.3)(6.0,5.2)(4.0,5.6)(1.3,5.0)(0.2,3.4)(1.2,1.9)(3.2,1.35)
\Linse
\rput(7.2,0){\Linse}
    \rput(14.3,1.8){$M$}
    \pscurve[linewidth=2pt](3.6,1.35)(3.4,2.4)(3.6,3.05)
    \pscurve[linewidth=2pt](11.2,1.35)(11.0,2.4)(11.2,3.05)
    \rput(3.6,1){$\gamma_1$}
    \rput(11.2,1){$\gamma_2$}
% tilde M oben
\rput(0,9){
    \pscurve(3.2,1.35)(5,1.5)(6.2,2.1)(7.4,2.8)(8.6,2.1)(9.8,1.5)(11.6,1.35)%
    (13.6,1.9)(14.6,3.4)(13.5,5.0)(10.8,5.6)(8.8,5.2)(7.4,4.3)(6.0,5.2)(4.0,5.6)(1.3,5.0)(0.2,3.4)(1.2,1.9)(3.2,1.35)
    \Linse    
    \rput(7.2,0){\Linse}
    \rput(14.3,1.8){$\widetilde M$}
    \psframe*[linecolor=white](3.5,1.2)(3.7,1.4)
    \psframe*[linecolor=white](3.5,2.9)(3.7,3.2)
    \pscurve(3.5,1.35)(3.3,2.4)(3.5,3.05)
    \pscurve(3.7,1.35)(3.5,2.4)(3.7,3.05)
    \psframe*[linecolor=white](3.5,1.2)(3.7,1.4)
    \psframe*[linecolor=white](3.5,2.9)(3.7,3.2)
    \pscurve(3.5,1.35)(3.3,2.4)(3.5,3.05)
    \pscurve(3.7,1.35)(3.5,2.4)(3.7,3.05)
    \psframe*[linecolor=white](11.1,1.2)(11.3,1.4)
    \psframe*[linecolor=white](11.1,2.9)(11.3,3.2)
    \pscurve(11.1,1.35)(10.9,2.4)(11.1,3.05)
    \pscurve(11.3,1.35)(11.1,2.4)(11.3,3.05)
   }
% Der Pfeil
\psline[linewidth=1.5pt]{->}(7.4,9.5)(7.4,7)

% Kurve, die rumgeht
%    \pscurve[linewidth=2pt](3.5,2.8)(4.6,2.8)(6.0,3.9)(3.8,5.0)(2.0,4.58)(1.1,3.4)(2,2.4)(3.0,2.2)(3.4,2.2)
\end{pspicture}
\end{center}
\caption{The cut-open $\witi M$ and its projection onto $M$}
\end{figure}

%%%%
%
%%%%

Riemannian metrics and spin structures on $M$ can be pulled back to $\witi M$.

\begin{lemma}\label{cutopenlemma}
Let $\ga_1,\dots,\ga_g$ be a cut of $M$. 
Then the cut-open $\witi M$ is diffeomorphic to a sphere $S^2$ with $2g$ 
disks removed.
Moreover, if it is a spin-cut, $\witi M$ carries the spin structure inherited 
from $S^2$.
\end{lemma}

\proof{}
At first we prove that $\witi M$ is connected. 
Assume that $\witi M$ is not connected. This would imply that the boundary of one
of the connected components of $\witi M$ is homologous to zero. Hence a nontrivial linear combination
of the $[\ga_i]$ vanishes which is impossible by the definition of a cut.

Since the Euler characteristic of $\witi M$ satisfies 
 $\chi(\witi M)=\chi(M)=2-2g$ and $\witi M$ has $2g$ boundary circles, it must
be diffeomorphic to a sphere  $S^2$ with $2g$ disks removed.

In the case of a spin-cut, the spin structure is nontrivial along each of the boundary components. Therefore the spin structure extends to the disk
which has been removed.
Hence $\witi M$ carries the spin structure which is the pullback of the unique spin structure on $S^2$ under any injective immersion 
$\witi M\hookrightarrow S^2$. 
\qed
 
\begin{definition}
Let $M$ be a closed surface with a fixed Riemannian metric and a fixed spin 
structure with Arf invariant $1$.
Let $\ga_1,\dots,\ga_g$ be a spin-cut. 
Denote by $\pa_1\witi M,\dots, \pa_{2g}\witi M$ the boundary components
of the cut-open $\witi M$.
We define the \emph{cut-diameter} of the spin-cut by
  $$\de(\ga_1,\dots,\ga_g):=\min_{1\leq i< j\leq 2g} d\big(\pa_i\witi M,\pa_j\witi M\big),$$
where $d(A,B)$ denotes the length of a shortest path joining $A$ and $B$.
The \emph{ spin-cut-diameter} of $M$ is defined as 
  $$\de(M):=\sup \de(\ga_1,\dots,\ga_g)$$
with the supremum running over all spin-cuts. 
The  spin-cut-diameter $\de(M)$ is a finite positive number
depending on the surface $M$, the Riemannian metric and the spin structure.
\end{definition}

%%%%%%%%%%%
% 
% picture: cutdiam.pic
% 
%%%%%%%%%%%

\begin{figure}[h]
\begin{center}
\psset{unit=7mm}
\begin{pspicture}(1,1)(14,6)
%   \psset{showpoints=true}
%    \psgrid(-2,-2)(18,10)

%% M unten
    \psccurve(3.2,1.35)(7.4,2.4)(11.6,1.35)%
    (14.6,3.4)(11.8,5.6)(7.4,4.9)(3.0,5.6)(0.2,3.4)
    \pscurve(1.8,4.0)(2.0,3.7)(3.7,3.05)(5.4,3.7)(5.6,4.0)
    \pscurve(2.0,3.7)(3,4.2)(4,4.25)(5.1,4.0)(5.4,3.7)
\rput(7.2,0)
      {\pscurve(1.8,4.0)(2.0,3.7)(3.7,3.05)(5.4,3.7)(5.6,4.0)
      \pscurve(2.0,3.7)(3,4.2)(4,4.25)(5.1,4.0)(5.4,3.7)}
    \pscurve[linestyle=dotted,linestyle=dotted](3.52,1.7)(4.4,1.75)(7.4,2.8)(10.4,1.75)(11.28,1.7)
    \psecurve[linestyle=dotted](2.6,3.0)(3.5,2.7)(6.1,3.8)(3.6,4.9)(1.4,3.8)(3.4,2.2)(4.0,2.5)
\rput(7.2,0)
      {\psecurve[linestyle=dotted](2.6,3.0)(4.15,2.7)(6.1,3.8)(3.6,4.9)(1.4,3.8)(4.2,2.2)(4.5,2.5)
      }
    \psecurve[linestyle=dotted](11,1.95)(11.3,1.9)(13.8,3.8)(10.8,5.2)(8.6,4.6)(7.4,3.6)(4.4,2.4)(3.4,2.4)(3.1,2.5)
    \rput(14.3,1.8){$M$}
    \pscurve[linewidth=1.5pt](3.6,1.35)(3.4,2.4)(3.6,3.05)
    \pscurve[linewidth=1.5pt](11.2,1.35)(11.4,2.4)(11.2,3.05)
    \rput(3.6,1){$\gamma_1$}
    \rput(11.2,1){$\gamma_2$}

\end{pspicture}
\end{center}
\caption{The cut-diameter is the length of the shortest dotted line
(only representatives of 4 of the 6 homotopy classes of lines are shown)}
\end{figure}
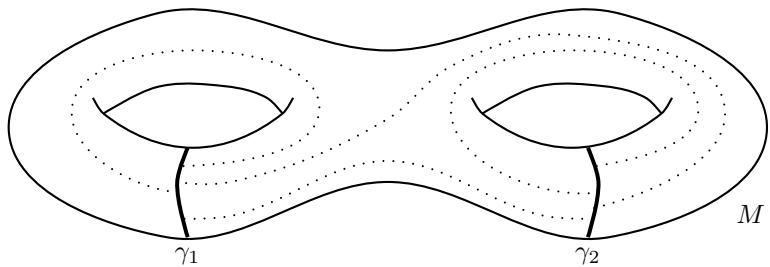

%%%%%%%%%%%%%%%%%%%%%%%%%%%%%%%%%%%
\section{Stable norms and the spin-cut-diameter}
\label{stable}
%%%%%%%%%%%%%%%%%%%%%%%%%%%%%%%%%%%

Let $M$ be a closed Riemannian manifold. 
In this section we define norms on $H_1(M,\RR)$ and $H^1(M,\RR)$, 
the stable norms, and we recall some of their properties. 
We will be able to express the  spin-cut-diameter defined in the
previous section in terms of stable norms of certain cohomology classes
which depend on the spin structure.
A good reference for stable norms is \cite{gromov:99}, Chapter 4C.
A more detailed exposition of stable norms can be found in 
\cite{federer:96}. 

For any $v\in H_1(M,\RR)$ the \emph{stable norm} is defined as
  $$\stabnorm{v}:=\inf \left\{ \sum_{i=1}^k |a_i|\cdot\length{c_i}\right\}$$
where the infimum runs over all 1-cycles
$\sum_{i=1}^k a_i c_i$ representing $v$ with $a_i\in \RR$, 
$k\in \NN\cup\{0\}$ and $c_i:S^1\to M$ smooth. 

For cohomology classes $\al\in H^1(M,\RR)$ we define the \emph{stable norm} by
  $$\stabnorm{\al}:=\inf \|\om\|_{L^\infty},$$
where the infimum runs over all closed smooth 1-forms $\om$ representing $\al$.

These norms are dual to each other in the following sense:
  $$\stabnorm{\al}
     = \sup \left\{\al(v)\, | \, v \in H_1(M,\RR), \stabnorm{v}=1\right\},$$
  $$\stabnorm{v}
     =  \sup \left\{\al(v)\, | \, \al \in H^1(M,\RR), \stabnorm{\al}=1\right\}.$$

We can also characterize the stable norm on $H_1(M,\RR)$ in terms of lengths 
of closed curves. For any 1-cycle $v\in H^1(M,\RR)$ which lies in the image of the map
$H^1(M,\ZZ)\to H^1(M,\RR)$
the relation
\begin{eqnarray*}
  \stabnorm{v} & = &\inf \Big\{{1\over n}\,\length{\ga}\,\Big|\, \ga
    \mbox{ is a closed curve representing }nv, n\in \NN\Big\}
\end{eqnarray*}
holds.

If $M=T^n$, the $n$-dimensional torus with an arbitrary Riemannian metric, 
then one can identify $H_1(T^n,\RR)$ with the universal covering of $T^n$.
Let $d$ be the distance function on $H_1(T^n,\RR)$ induced by the pullback
of the Riemannian metric on $T^n$.
Burago \cite{burago:92} proved that there is a 
constant $C$, such that for any $x,y\in H_1(T^n,\RR)$
  $$|d(x,y)-\stabnorm{x-y}|\leq C.$$

Roughly speaking, this result says that the stable norm is a good 
approximation for the distance $d$.

The stable norm also plays a central role in Bangert's criterion 
\cite{bangert:90} for the existence of globally minimizing 
geodesics on the universal covering 
$\widetilde{M}$ of a closed Riemannian manifold $M$. 
E.~g.\ if $b_1(M)\geq 2$, and if the stable norm on $H_1(M,\RR)$ is strongly 
convex, then there are infinitely many geodesics on $M$ whose lifts
are globally minimizing geodesics on $\widetilde{M}$.   

In the special case that $M$ is a closed orientable surface of positive genus,
any closed curve $\ga$ representing a nontrivial 
$[\ga]=[\al]^n\in \pi_1(M)$ with $n\geq 2$
has a self-intersection. To see this, let $\overline{M}$ be the universal 
covering. We lift $\ga$ to $\overline{M}/\lan [\al]\ran$ where $[\al]$
acts via deck transformations and apply 
Lemma~\ref{appendixlemmaA} for $S\setminus\{N,S\}\cong\overline{M}/\lan [\al]\ran$. 
A standard curve shortening argument
shows that in this case we can characterize the stable norm 
of an integral class $v$ as follows:
\begin{eqnarray*}
\stabnorm{v} &=& 
\inf\big\{\length{\ga}\,\big|\, 
\mbox{$\ga$ is a closed curve in $M$ representing $v$}\big\} .
\end{eqnarray*}

\begin{remark}
An intersection argument implies that $\stabnorm{\cdot}$
is a strictly convex norm on $H_1(T^2,\RR)$ \cite{massart:97a}.
In contrast to this, on any surface of genus $\geq2$ the stable norm is 
not strictly convex \cite{massart:97a}.
\end{remark}

\komment{
For $M=T^2$ the map $i:H^1(T^2,\ZZ)\to H^1(T^2,\RR)$ is an injective,
whereas $p:H^1(T^2,\ZZ) \to H^1(T^2,\ZZ_2)$.
For $\al\in H^1(T^2,\ZZ_2)$ the definition
  $$\stabnorm{\al}:=\left\{\stabnorm{i(\hat\al)}\,\Big|\,\hat\al\in H^1(T^2,\ZZ),
    j(\hat\al)=\al\left\}$$
yields a norm on  $H^1(T^2,\ZZ_2)$, the \emph{stable norm on $H^1(T^2,\ZZ_2)$}.
}

In the remaining part of this section we specialize to the case $M=T^2$, and 
we will show how the stable norm can be used to 
express the  spin-cut-diameter of a spin structure.

Let $\ga:S^1\to T^2$ be a noncontractible simple closed curve along which the spin structure
is nontrivial. Then $[\ga]\in H_1(T^2,\ZZ)\setminus \{0\}$.
We define $\al_\ga\in H^1(T^2,\ZZ)$ via the relation
  $$\lan\al_\ga,\be\ran=[\ga]\cap \be, \qquad \forall \be \in H_1(T^2,\ZZ).$$

\begin{proposition}\label{dusys}
Let $\de(M)$ be the  spin-cut-diameter of a 2-torus with spin structure
whose Arf invariant equals $1$.  
Let $\ga_0:S^1\to T^2$ be a noncontractible simple closed curve 
along which the spin structure is nontrivial, \ie $\ga_0$ is a spin-cut of $M$.
Then for
\begin{eqnarray*}
\de_0:=\sup\{\de(\ga)\,|\,\ga \mbox{\rm\ is a simple closed curve homotopic to }\ga_0\}
\end{eqnarray*}
we have 
  $$\de_0=\frac{1}{\stabnorm{\al_{\ga_0}}}.$$
\end{proposition}

\proof{}
\begin{enumerate}[(a)]
\item\label{kleiner} We show $\de_0\leq1/\stabnorm{\al_{\ga_0}}$.

Let $\ep>0$.
Choose a simple closed curve $\ga$ homotopic to $\ga_0$ 
such that $\de(\ga)\geq (1+\ep)^{-1}\de_0$.
We cut $T^2$ along $\ga$. 
Then the cut-open $\ccutM$ thus obtained is a topological cylinder.
Let $\ti c:[a,b]\to \ccutM$ be a curve of minimal length joining the two boundary components 
$\pa_1\ccutM$ and $\pa_2\ccutM$ of $\ccutM$. 
Let $c$ be the image of $\ti c$ under the map $\ccutM\to T^2$. 
Clearly $\length{c}=\length{\ti c}=\de(\ga)\geq (1+\ep)^{-1}\de_0$.
Let $f:\ccutM\to [0,\de_0]$ be a smooth function 
with the following properties:
  $$\matrix{|df|\leq 1+2\ep,\hfill\cr
      \matrix{
          f\equiv 0 \hfill&\mbox{on a neighborhood of } \pa_1\ccutM,\hfill\cr
          f\equiv \de_0\hfill&\mbox{on a neighborhood of } \pa_2\ccutM.\hfill
      }
    } 
  $$  
Such an $f$ can be obtained for example by a smooth approximation
of the Lipschitz function 
\begin{eqnarray*}
  \bar{f}:\ccutM&\to& [0,\de_0],\\
  x&\mapsto& \,{\de_0\over\de(\ga)} \, \min\Big\{d(x,\pa_1\ccutM),\de(\ga)\Big\}.
\end{eqnarray*}

Let $\om$ be the $1$-form on $T^2$ such that $df$ equals the pullback of 
$\om$.  

We now prove $\de_0\cdot\al_\ga=\pm [\om]$.

Observe that $\om(\dot \ga(t))={d\over dt}\, \left(f \circ \ga\right) 
\equiv 0$, since $f$ is constant along $\pa_1\ccutM$.
Hence $\int_{\res{\ga}{I}}\om =0$ for any $I\subset S^1$. 
In particular, 
  $$\lan [\om],[\ga]\ran=\int_\ga \om =0.$$ 
There are $t_1,t_2\in S^1$ such that 
$\ga(t_1)= c(a)$, $\ga(t_2)= c(b)$. 
Let $\be$ be the product path $\be:= \res{\ga}{[t_2,t_1]}* c$. 

%%%%%%%%
%
% picture sytor.pic
% 
%%%%%%%%

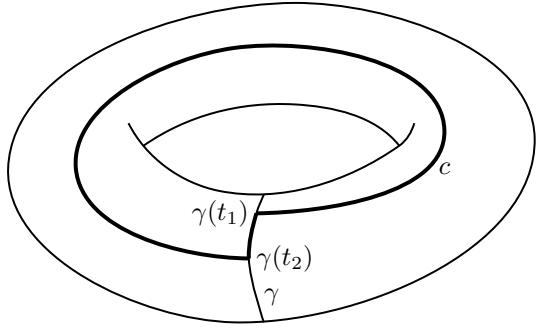
\begin{figure}[h]
\begin{center}
\begin{pspicture}(-0.4,1)(8,6.3)
\psset{linewidth=.75pt}
%    \psset{showpoints=true}    
%   \psgrid(-2,-2)(13,20)
    \psccurve(3.2,1.35)(7.2,3.8)(4.0,5.6)(0.2,3.4)
    \psecurve(1.8,5)(1.8,4.0)(3.0,3.1)(5.4,3.7)(5.6,4.0)(5.6,5)
    \psecurve(1.5,3)(2.0,3.7)(4,4.25)(5.4,3.7)(5.5,3)
    \pscurve(3.6,1.35)(3.4,2.2)(3.5,2.8)(3.6,3.05)
    \psecurve[linewidth=1.5pt](3.6,1.35)(3.4,2.2)(3.5,2.8)(3.6,3.05)
    \psecurve[linewidth=1.5pt](3,3)(3.5,2.8)(6.0,3.9)(3.2,5.0)(1.1,3.4)(3.4,2.2)(4,2.5)

%Beschriftungen
\rput[lt](3.6,1.8){$\ga$} 
\rput[l](3.5,2.2){$\ga(t_2)$} 
\rput[r](3.4,2.8){$\ga(t_1)$} 
\rput(6.0,3.4){$c$} 
\end{pspicture}
\end{center}
\caption{The curve $\be$ in the proof of Proposition~\ref{dusys} (thick line)}
\end{figure}

%%%%%%%%%%%

Then $[\ga]\cap[\be]=\pm 1$.
Moreover, 
\begin{eqnarray*}
\lan [\om],[\be]\ran&=&\underbrace{\int_{\res{\ga}{[t_2,t_1]}}\om}_{=0} 
    + \int_{c}\om
\;=\;
\int_{c}df\;=\; f(c(b))-f(c(a))\\
&=&
\de_0 \;=\; \pm \de_0\, [\ga]\cap[\be] \;=\; \pm \de_0\, \lan \al_\ga , [\be]\ran .
\end{eqnarray*}
Therefore $[\om]\mp \de_0\cdot\al_\ga$ vanishes on $[\ga]$ and on $[\be]$.
Since $[\ga]$ and $[\be]$ form a basis of $H_1(T^2,\ZZ)$ we obtain 
$\de_0\cdot\al_\ga=\pm [\om]$.

From 
$$
\de_0\cdot\stabnorm{\al_\ga} = \stabnorm{[\om]}\le\|\om\|_{L^\infty} 
\leq 1+2\ep
$$
we get the $\leq$-part of the equation by taking the limit $\ep\to 0$.

\item
Now we prove $\de_0\geq1/\stabnorm{\al_{\ga_0}}$.

We choose a smooth closed 1-form $\om$ on $T^2$ such that $[\om]=\al_{\ga_0}$
and $\|\om\|_{L^\infty}\leq \stabnorm{\al_{\ga_0}}+\ep$ for small $\ep>0$.
The cyclic subgroup $\lan[\ga_0]\ran$ of $\pi_1(T^2)$ generated by $[\ga_0]$
acts via deck transformations on the universal covering $\RR^2$ of $T^2$. 
Define the cylinder $Z:=\RR^2/\lan[\ga_0]\ran$. 
Since $[\ga_0]$ generates the first cohomology of $Z$ and $\al_{\ga_0}$
vanishes on $[\ga_0]$ the pullback of the cohomology class $[\om]=\al_{\ga_0}$
is trivial on $Z$.
Hence we can find a smooth function $f:Z\to \RR$ such that $df$ is the pullback
of $\om$ under the covering $Z\to T^2$. 

The function $f$ is proper.
Without loss of generality we can
assume that~$0$ is a regular value of $f$. 
Then $f^{-1}(0)$ is a union of simple closed curves.
According to Lemma~\ref{appendixlemmaB} there is a simple closed curve
$\ga$ in $f^{-1}(0)$ whose homotopy class
generates $\pi_1(Z)$. Choose the orientation of $\ga$ such that  
$\ga$ is homotopic to $\ga_0$.   
The spin struture is nontrivial along $\ga$, 
hence $\ga$ defines a spin-cut $\witi M \to M$, \ie a map which 
is a diffeomorphism from the interior of $\witi M$ onto $M\setminus\ga(S^1)$
and a trivial double covering from $\pa \witi M$ onto $\ga(S^1)$.

We can identify $\witi M$ with a closed subset of $Z$, and we can assume 
that $\res{f}{\pa_1\ccutM}\equiv 0$, $\res{f}{\pa_2\ccutM}\equiv 1$,
where $\pa_1\ccutM$ and $\pa_2\ccutM$ denote the two boundary components 
of $\witi M$.

Let $c:[a,b]\to \witi M$ 
be a curve of minimal length joining the two boundary components 
$\pa_1\ccutM$ and $\pa_2\ccutM$. By definition we have
$\de(\ga)=\length{c}$. It follows
\begin{eqnarray*}
  1 &=& f(c(b))-f(c(a)) = \int_{c} \,df\\ 
    &\leq &\length{c} \,\|df\|_{L^\infty}\\
    & = & \de(\ga)\,\|\om\|_{L^\infty}\\
    &\leq & \de_0 \, \Bigl(\stabnorm{\al} + \ep\Bigr).
\end{eqnarray*}

The limit $\ep\to 0$ yields $\de_0\ge/\stabnorm{\al_{\ga_0}}$.
\end{enumerate}
\vskip-5mm\qed

\begin{corollary}
The  spin-cut-diameter satisfies
\begin{eqnarray*}
\de(M)&=&\sup\, \Biggl\{{1\over \stabnorm{\al_\ga}}\;\Big|\;\ga \mbox{\/\rm\  is a noncontractible
simple closed curve} \\
&&\phantom{\sup\,\Biggl\{}\mbox{\rm along which the spin structure is nontrivial.}\Biggr\}
\end{eqnarray*}
\qed
\end{corollary}

%%%%%%%%%%%%%%%%%%%%%%%%%%%%%%%%%%%
\section{An estimate for the 2-torus}
\label{torussec}
%%%%%%%%%%%%%%%%%%%%%%%%%%%%%%%%%%%

We now come to the first main result of this paper.
We give a geometric lower bound for the eigenvalues of 
the Dirac operator on a 2-torus which is nontrivial for all metrics
and for all spin structures.

\begin{theorem}\label{2torus}
Let $T^2$ be the 2-torus equipped with an arbitrary Riemannian metric and 
a spin structure whose Arf invariant equals $1$. 
Let $\la$ be an eigenvalue of the Dirac operator and let $\de(T^2)$ be
the spin-cut-diameter. 
Then for any $k\in \NN$, 
$$
|\la|\geq -{2\over k\,\de(T^2)} + \sqrt{{\pi\over k\,\area(T^2)}+{2\over k^2 \de(T^2)^2}}.
$$
\end{theorem}

Note that the right hand side of this inequality is positive for sufficiently 
large~$k$, but tends to $0$ for $k\to \infty$.
The best bound is obtained by choosing
  $$k=\left[4\,(1+\sqrt{2})\,{\area(T^2)\over \pi\,\de(T^2)^2}\right]$$
or 
  $$k=\left[4\,(1+\sqrt{2})\,{\area(T^2)\over \pi\,\de(T^2)^2}\right]+1.$$

\proof{}
Let $\ga$ be a spin-cut, \ie $\ga$ is a simple closed curve in $T^2$ along 
which the spin structure is nontrivial.
Assume $\de(\ga)\geq (1+\ep)^{-1}\,\de(T^2)$ for small $\ep>0$.
 
We now proceed as in part 
(\ref{kleiner}) of the proof of Proposition~\ref{dusys}.
On the cut-open $\witi T^2$ we obtain a function $f:\witi T^2\to [0,\de(T^2)]$
satisfying
  $$\matrix{|df|\leq 1+2\ep,\hfill\cr
      \matrix{
          f\equiv 0 \hfill&\mbox{on a neighborhood of } \pa_1\witi T^2,\hfill\cr
          f\equiv \de(T^2)\hfill&\mbox{on a neighborhood of } \pa_2\witi T^2.\hfill
      }
    } 
  $$ 
Let $\om$ be the $1$-form on $T^2$ such that $df$ equals the pullback of $\om$.  
 
The homotopy class $[\ga]\in \pi_1(T^2)$ acts on the universal covering 
$\RR^2$ of $T^2$, and
  $$Z:=\RR^2/ \lan[\ga]\ran$$
is a cylinder covering $T^2$.
We pull the metric and the spin structure on $T^2$ back to a metric and 
a spin structure on $Z$.
 
We fix a $w\in\pi_1(T^2)$ with $[\ga]\cap w=1$. Then $w$
generates the deck transformation group of the covering $Z\to T^2$. 
Let $\ti\ga:S^1\to Z$ be a lift of $\ga$.   
Then $Z\setminus \big(\ti\ga(S^1)\cup w\cdot \ti\ga(S^1)\big)$ consists
of three connected components. 
Two of them are unbounded and one is bounded. 
The closure of the bounded component can be identified 
with the cut-open $\witi T^2$. 
The function $f$ can then be extended ``pseudo-periodically'' to $Z$, 
more precisely, 
\begin{equation}\label{period}
  f(w+p)=\de(T^2) +f(p)
\end{equation}
for all $p\in Z$, where $w$ acts as a deck transformation on $Z$.
Note that 
  $$\area\Big(f^{-1}\big((t,t+\de(T^2)]\big)\Big)=\area\big(\witi T^2\big)=\area(T^2).$$

%%%%%%%%%%%%%
%
% picture cylcovering.pic
% 
%%%%%%%%%%%%%

\begin{figure}[h]
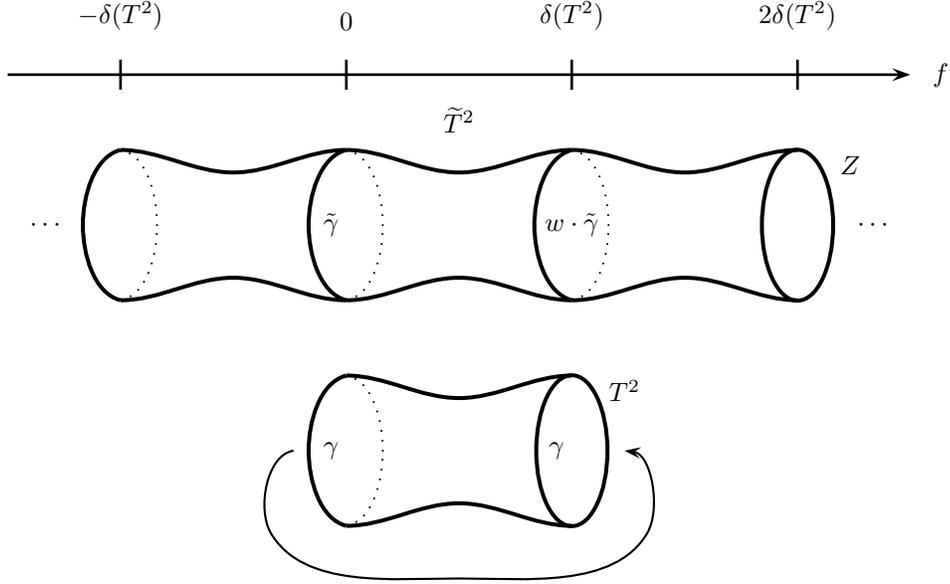

\begin{center}
\pspicture(-6,-1)(6,7)
\psset{unit=1cm}
%\psgrid(-6,-1)(6,7)

\def\roehre{
\psellipse[linestyle=dotted](-1.5,1)(0.5,1)
\psecurve[linewidth=1.5pt](-3,0.3)(-1.5,0)(0,0.3)(1.5,0)(3,0.3)
\psecurve[linewidth=1.5pt](-3,1.7)(-1.5,2)(0,1.7)(1.5,2)(3,1.7)
\psecurve[linewidth=1.5pt](1.5,0)(-1.5,2)(-2,1)(-1.5,0)(1.5,2)}

\rput(0,0){\roehre}
\psellipse[linewidth=1.5pt](1.5,1)(0.5,1.03)
\psecurve{->}(-1,-1)(-2.2,1)(-2.5,-0.1)(0,-0.7)(2.5,-0.1)(2.2,1)(1,-1)
\rput(2.2,1.8){$T^2$}

\rput(0,3){\roehre}
\rput(-3,3){\roehre}
\rput(3,3){\roehre}
\psellipse[linewidth=1.5pt](4.5,4)(0.5,1.03)
\rput(5.5,4){$\cdots$}
\rput(-5.5,4){$\cdots$}
\rput(5.2,4.8){$Z$}

\rput(0,5.4){$\witi T^2$}
\rput(-1.7,4){$\ti \ga$}
\rput(-1.7,1){$\ga$}
\rput(1.3,1){$\ga$}
\rput(1.5,4){$w\cdot\ti \ga$}
\rput(0,6){
  \psline[linewidth=1pt]{->}(-6,0)(6,0)
  \rput[l](6.3,0){$f$}
  \def\strich#1{\psline[linewidth=1pt](0,-.2)(0,.2)
              \rput[b](0,.6){$#1$}
             }
  \rput(-4.5,0){\strich{-\de(T^2)}}
  \rput(-1.5,0){\strich{0}}
  \rput(1.5,0){\strich{\de(T^2)}}
  \rput(4.5,0){\strich{2\de(T^2)}}
}

\endpspicture
\end{center}
\caption{The cylinder $Z$ and a fundamental domain}
\end{figure}

%%%%%%%%%%%%%%%%%%%%%

We set 
    $$\fm:=f^{-1}\Big([-k\de(T^2),0]\Big),$$ 
    $$\fp:=f^{-1}\Big([0,k\de(T^2)]\Big).$$ 
Both $\fm$ and $\fp$ are isometric to $k$ copies of $\witi T^2$ glued together
to a cylinder.
Similarly, we consider $\fm\cup \fp$ as a cylinder consisting of $2k$ copies 
of $\witi T^2$.
We glue two disks to the remaining two boundary components of $\fm\cup \fp$
and obtain a surface $N$ of genus $0$. 
We extend the metric on $\fm\cup \fp$ to one on $N$ such that the total area 
of the two disk glued in is smaller than $\ep$.
Hence
$$
\area(N)\leq 2k\, \area(T^2)+\ep.
$$

By Proposition~\ref{nontrivcrit} the spin structure on 
$\fm\cup \fp$ extends to the unique spin structure on $N$.

For fixed $k\in \NN$ let $X_1:\RR\to[0,1]$ be a smooth function with
\begin{eqnarray*}
   X_1(t)=1 & \mbox{for} & t\leq 0,\\
   X_1(t)=0 & \mbox{for} & t\geq k,
\end{eqnarray*}
$$
   |X_1'(t)|\leq {1+\ep\over k}  \mbox{\ \ for all $t$.\ \ }
$$
%  $$X(t):=\cases{1-|t|/k& for $-k\leq t\leq 0$\cr                  
%                 1-t/k & for $0\leq t\leq k$\cr
%                 0& elsewhere} $$

%%%%%%%%%%%%%%%%%%%%%%%%%%%%%%%%%%
% Funktionsgraph
%%%%%%%%%%%%%%%%%%%%%%%%%%%%%%%%%%
%%%
% 
% picture: graph.pic
%
%%%%%%%%%%%%%%%%%%%%%%%%%%%%%%%%%%

\begin{figure}[h]
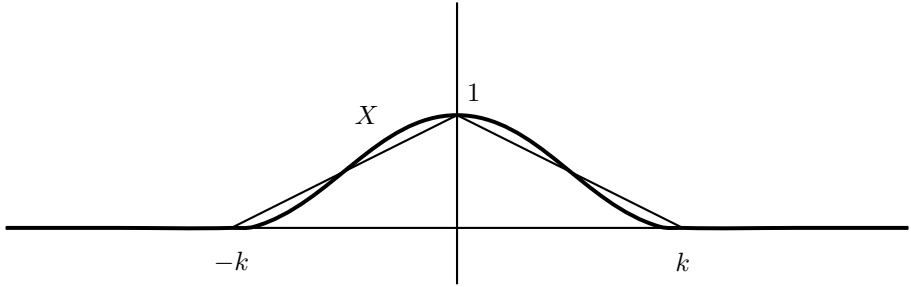

\begin{center}
\pspicture(-6,-1)(6,4)
\psset{unit=1.5cm}
%\psgrid(-6,-1)(6,4)

\psline(-4,0)(4,0)
\psline(0,-0.5)(0,2)
\psline(-2,0)(0,1)
\psline(2,0)(0,1)
\psecurve[linewidth=1.5pt](-5,0)(-4,0)(-3,0)(-2,0)(-1.8,0.01)(0,1)(1.8,0.01)(2,0)(3,0)(4,0)(5,0)

\rput(-2,-0.3){$-k$}
\rput(2,-0.3){$k$}
\rput(0.15,1.2){$1$}
\rput(-0.8,1){$X$}

\endpspicture
\end{center}
\caption{The graph of $t\mapsto X(t)$}
\end{figure}

%%%%%%%%%%%%%%%%%

We set $X(t):=X_1(t)-X_1(t+k)$.
Then 
  $$\chi(p):=X\left({ f(p)\over \de(T^2)}\right)$$ 
is a compactly supported smooth function on $Z$ 
with 
$$
k\cdot \|\na \chi\|_{L^\infty} \le
k\cdot \| X'\|_{L^\infty}\cdot \frac{\| df\|_{L^\infty}}{\de(T^2)}
\leq {(1+\ep)(1+2\ep)\over \de(T^2)}=:a_\ep.
$$
We denote the $L^2$-norm of a spinor $\ph$ on a subset $A$ of
the manifold on which $\ph$ is defined by
$$
\|\ph\|_A := \sqrt{\int_A |\ph|^2 \;\mbox{d\,area}} .
$$
If $A$ equals the whole manifold we simply write
$$
\|\ph\|_A =: \|\ph\| .
$$
Now let $\ph$ be an eigenspinor on $T^2$ corresponding to an eigenvalue $\la$ 
of the Dirac operator.
By the preceeding lemma, 
the spin structure pulled back via $\pi$ extends to the unique spin structure on $N$. 
Thus $\chi\cdot\pi^*\ph$ is a well-defined spinor on $N$, and 
we obtain the following estimate
\begin{eqnarray}
   \|D(\chi\cdot\pi^*\ph)\|_{\fm}^2 
     & = & \|\na \chi\cdot\pi^*\ph+ \chi\cdot D(\pi^*\ph)\|_{\fm}^2\nonumber\\
     & \leq & \left({a_\ep\over k}\cdot \|\pi^*\ph\|_{\fm} + |\la|\,\|\chi\cdot \pi^*\ph\|_{\fm}\right)^2\nonumber\\
     & \leq & {a_\ep^2\over k^2}\cdot \|\pi^*\ph\|_{\fm}^2 
     + {2|\la|a_\ep \over k}\|\pi^*\ph\|^2_{\fm}\nonumber\\
     && 
     + \la^2\,\|\chi\cdot \pi^*\ph\|_{\fm}^2\nonumber\\
     & = & \left({a_\ep^2\over k} + 2|\la|\,a_\ep\right) \|\ph\|^2_{T^2}
     + \la^2\,\|\chi\cdot \pi^*\ph\|_{\fm}^2.\label{e1}
\end{eqnarray}
In a similar manner we obtain
\begin{eqnarray}
   \|D(\chi\cdot\pi^*\ph)\|_{\fp}^2 
     & \leq & \left({a_\ep^2\over k} + 2|\la|\,a_\ep\right) \|\ph\|^2_{T^2}
     + \la^2\,\|\chi\cdot \pi^*\ph\|_{\fp}^2.\label{e2}
\end{eqnarray}

From
  $$X(t)^2+X(t-k)^2=X(t)^2+ (1-X(t))^2\in [1/2,1]$$ 
for $0\leq t\leq k$ 
we obtain
\begin{eqnarray*}
  {k\over 2} \,\|\ph\|^2_{T^2}&\leq& \|\chi\cdot\pi^*\ph\|_{\fmp}^2 
      \leq  k \|\ph\|^2_{T^2}
\end{eqnarray*}
which together with \eref{e1} and \eref{e2} gives
\begin{eqnarray*}
   \|D(\chi\cdot\pi^*\ph)\|_{\fmp}^2 
     & \leq & \left\{2\left({a_\ep^2\over k} + 2|\la|\,a_\ep\right) +k\cdot \la^2\right\}   \|\ph\|^2_{T^2}.
\end{eqnarray*}

We plug $\chi\ph$ into the Rayleigh quotient and use Theorem~\ref{baersphere} 
to get
\begin{eqnarray*}
{4\pi\over 2k\,\area(T^2)+\ep} & \leq & {4\pi\over \area(N)}\\
&\leq & {\|D(\chi\cdot\pi^*\ph)\|_{\fmp}^2 
\over \|\chi\cdot\pi^*\ph\|_{\fmp}^2 }\\
&\leq & {{2a_\ep^2/ k} + 4|\la|\,a_\ep+k\cdot \la^2\over k/2 }.
\end{eqnarray*}
Thus
\begin{eqnarray*}
{\pi\over2k\, \area(T^2)+\ep}&\leq & {a_\ep^2\over k^2} + {2|\la|\,a_\ep\over k}+{\la^2\over 2}.
\end{eqnarray*}
In the limit as $\ep\to 0$ we obtain
\begin{eqnarray*}
{\pi\over k\, \area(T^2)}&\leq & {2\over k^2\de(T^2)^2} + {4|\la|\over k\,\de(T^2)}+\la^2.
\end{eqnarray*}

Solving this inequality proves the theorem. \qed

%%%%%%%%%%%%%%%%%%%%%%%%%%%%%%%%%%%%%%%%%%%
\section{Compact Surfaces of higher genus}
\label{hochsec}
%%%%%%%%%%%%%%%%%%%%%%%%%%%%%%%%%%%%%%%%%%%

Using a similar technique we can also obtain a lower bound for the 
Dirac spectrum on closed surfaces $M$ of higher genus.

\begin{theorem}\label{ghoch}
Let $M$ be a closed surface of genus $g\ge 1$ with a Riemannian metric 
and a spin structure whose Arf invariant equals 1.
Let $\de(M)$ be the  spin-cut-diameter of $M$.
Then for all eigenvalues $\la$ of the Dirac operator we have
  $$|\la| \geq {2\sqrt{\pi}\over (2g+1)\,\sqrt{\area(M)}} -{1\over \de(M)}.$$
\end{theorem}

Note that on any closed oriented surface of genus $g\geq 1$ there is a 
Riemannian metric and a spin structure such that $\de(M)^2/\area(M)$ 
is arbitrarily large.
To see this take a suitable finite graph $\Ga$ embedded in $\RR^3$ and let $M$
be the boundary (smoothed out appropriately) of a tubular neighborhood of 
$\Ga$ of small tubular radius $r>0$.
Provide $M$ with the Riemannian metric and the spin structure induced
from $\RR^3$.
Then for $r\to 0$ the  spin-cut-diameter stays bounded while 
the area tends to $0$.
By Theorem~\ref{ghoch} the smallest eigenvalue of $D^2$ must then tend to 
$\infty$.
Hence any closed oriented surface carries a Riemannian metric and a spin 
structure such that the above estimate is not trivial.

%%%%%%%%%%%
%
% picture : umgebung.pic
% 
%%%%%%%%%%%

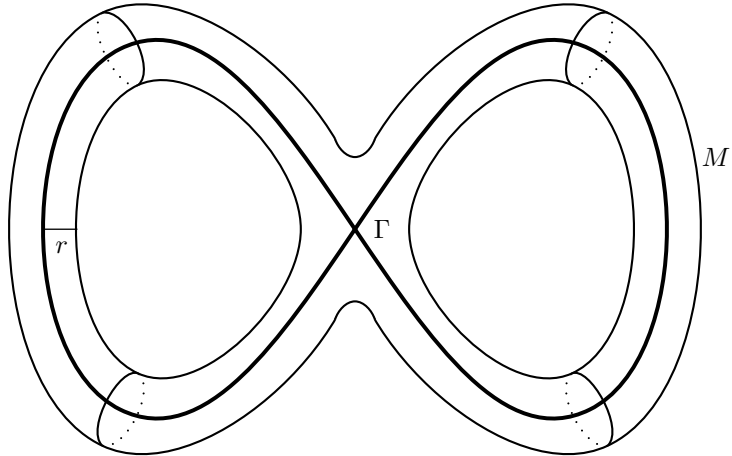
\begin{figure}[h]
\begin{center}
\psset{unit=12mm}
\begin{pspicture}(1,1.5)(11,7)
%  \psset{showpoints=true}
%  \psgrid(-2,-2)(18,10)

  \psccurve(3.2,1.6)(3.2,6.4)(6.1,4)
  \psccurve(8.8,6.4)(8.8,1.6)(5.9,4)
  \pspolygon[linecolor=white,fillstyle=solid,fillcolor=white](5.8,5)(6.2,5)(6.2,3)(5.8,3)
  \psecurve(6,6)(5.78,5)(6,4.8)(6.22,5)(6,6)
  \psecurve(6,2)(5.78,3)(6,3.2)(6.22,3)(6,2)
  \psccurve[linewidth=1.5pt](3.4,2)(3.4,6)(8.6,2)(8.6,6)
  \psccurve[linewidth=.8pt](5.4,4)(3.6,2.4)(3.6,5.6)
  \psccurve[linewidth=.8pt](6.6,4)(8.4,2.4)(8.4,5.6)

\psecurve(3,1)(3.6,2.4)(3.2,1.6)(4,2)
\psecurve[linestyle=dotted](3,2)(3.6,2.4)(3.2,1.6)(3,5)

\psecurve(9,1)(8.4,2.4)(8.8,1.6)(8,2)
\psecurve[linestyle=dotted](9,2)(8.4,2.4)(8.8,1.6)(9,5)

\psecurve[linestyle=dotted](3,7)(3.6,5.6)(3.2,6.4)(4,6)
\psecurve(3,6)(3.6,5.6)(3.2,6.4)(3,3)

\psecurve[linestyle=dotted](9,7)(8.4,5.6)(8.8,6.4)(8,6)
\psecurve(9,6)(8.4,5.6)(8.8,6.4)(9,3)

\psline[linewidth=0.5pt](2.56,4)(2.925,4)

\rput(10,4.8){$M$}
\rput(6.3,4){$\Gamma$}
\rput(2.75,3.8){$r$}

\end{pspicture}
\end{center}
\caption{The boundary $M$ of a small neighborhood of a graph $\Gamma$ 
in $\RR^3$ has a large 
maximal spin-cut-diameter compared to the area}
\end{figure}

%%%%%%%%%%%%%%%%%%%%

Theorem~\ref{ghoch} also holds for $g=1$ but in this case Theorem~\ref{2torus}
with $k=2$ gives a better estimate.

\proof{}
Let $\ga_1,\dots,\ga_g$ be a spin-cut of $M$.
We cut $M$ along the $\ga_i$ and obtain the cut-open $\witi M$.
According to Lemma~\ref{cutopenlemma}, $\witi M$ is a compact orientable surface
of genus~$0$ with $2g$ boundary components. 
The two boundary components of $\witi M$ that arise from cutting along 
$\ga_i$ we denote by $\pa_i^1\witi M$ and $\pa_i^2\witi M$.

We assume that the spin-cut has been chosen such that $\de(\ga_1,\dots,\ga_g)\geq \de(M)-\ep$ with $\ep>0$ small.
 
We take $2g+1$ copies of $\ccutM$, denoted by $\ccutM_0,\dots,\ccutM_{2g}$.
For $t=1,\dots,g$ we glue $\pa^1_t\ccutM_t$ to $\pa^2_t\ccutM_0$ and 
$\pa^2_t\ccutM_{g+t}$ to $\pa^1_t\ccutM_0$. The resulting surface $S_0$
is of genus~$0$ with $2g(2g-1)$ boundary components. We glue disks to these boundaries
and obtain a surface $S$ diffeomorphic to $S^2$.

%%%%%%%%
%
% picture ggetwo.pic
% 
%%%%%%%%

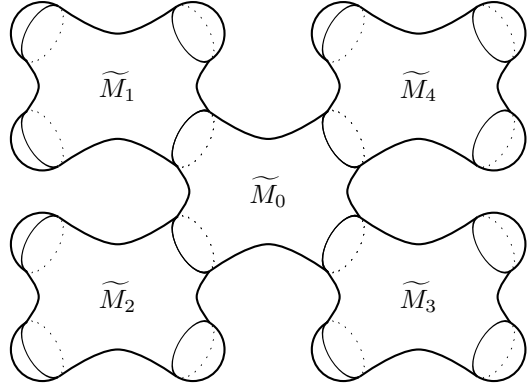
\begin{figure}[h]
\begin{center}
\psset{unit=10mm}
\begin{pspicture}(-1,-0.5)(7.5,5.5)
%   \psset{showpoints=true}
%   \psgrid(-1,-1)(10,10)

\def\xstueck#1#2
  {\rput{300}{
    \scalebox{1 .5}{
        \psarc[linestyle=#2](0,0){.42}{0}{180}
        \psarc[linestyle=solid](0,0){.42}{180}{360}
    }
  }
        \rput(1.0,.7){$\widetilde M_#1$}
\rput{240}(2.0,0){
    \scalebox{1 .5}{
        \psarc[linestyle=#2](0,0){.42}{0}{180}
        \psarc[linestyle=solid](0,0){.42}{180}{360}
    }
  }
\rput{240}(0,1.4){
    \scalebox{1 .5}{
        \psarc[linestyle=#2](0,0){.42}{0}{180}
        \psarc[linestyle=solid](0,0){.42}{180}{360}
    }
}
\rput{300}(2.0,1.4){
    \scalebox{1 .5}{
        \psarc[linestyle=#2](0,0){.42}{0}{180}
        \psarc[linestyle=solid](0,0){.42}{180}{360}
       }
    }
    \pscurve(.22,-0.36)(1,0)(1.78,-0.36)
    \pscurve(.22,1.76)(1,1.4)(1.78,1.76)
    \pscurve(-.22,.36)(-0.05,.7)(-.22,1.04)
    \pscurve(2.22,.36)(2.05,.7)(2.22,1.04)
}

\def\rokappe{%rechtsobenkappe
\psarc(2,1.4){.42}{300}{120}
}
\def\rukappe{
\psarc(2,0){.42}{240}{60}
}
\def\lokappe{
\psarc(0,1.4){.42}{60}{240}
}
\def\lukappe{%rechtsobenkappe
\psarc(0,0){.42}{120}{300}
}

\rput(0,0){\xstueck2{dotted} \lukappe\lokappe\rukappe}
\rput(2,1.4){\xstueck0{dotted} }
\rput(0,2.8){\xstueck1{dotted} \lukappe\lokappe\rokappe}
\rput(4,2.8){\xstueck4{dotted} \lokappe\rukappe\rokappe}
\rput(4,0){\xstueck3{dotted} \lukappe\rukappe\rokappe}
\end{pspicture}
\end{center}
\caption{The surface $S$ for $g=2$}
\end{figure}

%%%%%%%%%%%%%%%%%%%%%%%%

The Riemannian metric on $M$ pulls back to a Riemannian metric on 
$\ccutM$ and gives rise to a smooth metric on
$S_0$. 
We extend this metric to a metric on $S$ such that 
\begin{equation}
\area(S) \leq \area(S_0)+\ep =  (2g+1)\, \area(M)+\ep.
\label{eqq1}
\end{equation}
Since the spin structure of $M$ is nontrivial along each $\ga_i$, the induced
spin structures on $\ccutM_i$ fit together to the unique spin structure on $S$.

There is a smooth function $\chi\colon S\to [0,1]$ with the 
following properties:
\begin{enumerate}[(1)]
\item $\res{\chi}{\ccutM_0}\equiv 1$,
\item $\res{\chi}{S\setminus S_0}\equiv 0$,
\item $\|\na\chi\|_{L^\infty}\leq {\displaystyle{1\over \de(M)-2\ep}} $.
\end{enumerate}

Let $\ph$ be an eigenspinor of the Dirac operator on $M$ 
to the eigenvalue $\la$.
This spinor lifts to an eigenspinor $\ph_0$ of the Dirac operator on $S_0$.
Thus $\chi\cdot\ph_0$ is a well-defined spinor on $S$. We use it as a
test spinor 
for the Rayleigh quotient.
Theorem~\ref{baersphere} yields 
\begin{equation}\label{rqout}
  {4\pi\over \area(S)}\leq {\|D(\chi\cdot\ph_0)\|_S^2\over \|\chi\cdot\ph_0\|_S^2 .}
\end{equation}

We compute
$$
\|D(\chi\cdot\ph_0)\|_{\ccutM_i}^2
     \leq \left({1\over (\de(M)-2\ep)^2}+ |\la|\right)^2\|\ph\|_{M}^2 .
$$
Summing over $i$ yields
\begin{equation}
\|D(\chi\cdot\ph_0)\|_{S}^2\leq 
(2g+1)\,\left({1\over (\de(M)-2\ep)^2}+ |\la|\right)^2\|\ph\|_{M}^2.
\label{eqq2}
\end{equation}
The denominator of the Rayleigh quotient is estimated by
\begin{equation}
\|\chi\cdot\ph_0\|_S^2\geq \|\ph_0\|_{\ccutM_0}^2=\|\ph\|_M^2 .
\label{eqq3}
\end{equation}
Combining (\ref{eqq1}),(\ref{rqout}),(\ref{eqq2}), and (\ref{eqq3}) we obtain 
  $${4\pi\over  (2g+1)\,\area(M)+\ep}\leq (2g+1)\,\left({1\over (\de(M)-2\ep)^2}+ |\la|\right)^2$$
which yields 
in the limit $\ep\to 0$
  $${2\sqrt{\pi}\over (2g+1)\,\sqrt{\area(M)}} -{1\over \de(M)}\leq |\la|.$$
\qed

%%%%%%%%%%%%%%%%%%%%%%%%%%%%%%%%%%%
\section{An application to the Willmore integral}
\label{willmoresec}
%%%%%%%%%%%%%%%%%%%%%%%%%%%%%%%%%%%

The {\em Willmore integral} of an embedded closed surface $M\subset\RR^3$
is defined by
$$
W(M) = \int_M H^2 \mbox{dvol} = \| H \|^2
$$
where $H$ denotes the mean curvature of $M$.
The famous {\em Willmore conjecture} states that for an embedded 2-torus
the Willmore integral is bounded by
$$
W(M) \ge 2\pi^2 .
$$
This conjecture has been proven for various classes of embedded 2-tori
(see \cite{topping:00} for a good overview),
but in full generality it is still open.
We will not resolve this problem here but our estimates on Dirac eigenvalues
imply lower bounds on the Willmore integral as well.

Let $M\subset\RR^3$ be an embedded surface of genus $g \ge 1$.
The discussion from Sections~\ref{diracsurfsec} and \ref{spincutsection} shows
that the induced spin structure on $M$ admits spin-cuts and hence its
spin-cut-diameter $\de(M)$ is well-defined.
A spin-cut can be obtained by choosing disjoint simple closed curves 
$\ga_1,\ldots, \ga_g$ on $M$ which bound transversal disks in $\RR^3$
and whose homology classes $[\ga_1],\ldots, [\ga_g]$ in $H_1(M,\ZZ)$
are linearly independent.

\begin{theorem}
\label{willmore}
Let $T^2\subset\RR^3$ be an embedded torus.
Let $\de(T^2)$ be its spin-cut-diameter and let $W(T^2)$ be its Willmore
integral.
Then for any $k\in \NN$
$$
\sqrt{W(T^2)} \ge \sqrt{\frac{\pi}{k}+\frac{2\,\area(T^2)}{k^2\,\de(T^2)^2}}
- \frac{2\sqrt{\area(T^2)}}{k\,\de(T^2)}
$$
\end{theorem}

\proof{}
In \cite{baer:98} it was shown that a closed
surface possesses Dirac eigenvalues $\la$ satisfying
$$
\la^2 \le \frac{W(M)}{\area(M)}.
$$
Combining this with Theorem~\ref{2torus} yields the result.
\qed

This theorem yields a positive lower bound on $W(T^2)$ for all embedded
2-tori. 

\begin{remark}
From Theorem~\ref{ghoch} we can obtain a similar bound, but it turns out to 
be weaker than the well-known bound $W(M)\geq 4\pi$. 
\end{remark}

%%%%%%%%%%%%%%%%%%%%%%%%%%%%%%%%%%%
\section{Noncompact surfaces of finite area}
\label{noncompact}
%%%%%%%%%%%%%%%%%%%%%%%%%%%%%%%%%%%

Now we extend the bounds on Dirac eigenvalues to the $L^2$-spectrum
of the Dirac operator on a complete noncompact spin surface of finite area.
The {\em fundamental tone} of the square of the Dirac operator on a 
noncompact spin manifold is given by
$$
\la_\ast^2 = \inf_{\ph} \frac{\|D\ph\|^2}{\|\ph\|^2}
$$
where the infimum runs over all smooth spinors $\ph$ with compact support.
If $\la_\ast^2>0$, then the $L^2$-spectrum of $D$ has a gap about 0,
more precisely, 
$$
\lss(D) \cap (-\la_\ast,\la_\ast) = \emptyset.
$$
Any complete surface $M$ of finite area is diffeomorphic to a closed
surface $\barM$ with finitely many points removed.
The genus $g$ of $\barM$ is then also called the genus of~$M$.
By a {\em cut} of $M$ we mean a collection of simple closed curves
$\ga_1, \ldots ,\ga_g$ on $M$ which are mapped under the diffeomorphism
to a cut on $\barM$.
If $M$ carries a spin structure, then we call the cut a {\em spin-cut}
if the spin structure is nontrivial along all $\ga_i$ just as we did for
closed surfaces.
If the  spin structure on $M$ extends to one on $\barM$, then we say the
spin structure is {\em nontrivial along the ends}.

Given a spin-cut on $M$ one can define the {\em cut-open} as before.
It is now a noncompact complete surface of finite area with compact boundary.
The {\em spin-cut-diameter} is again defined as the minimal distance of
the various boundary components of the spin-cut.
Taking the supremum over all spin-cuts yields the 
\emph{spin-cut-diameter} $\de(M)$ depending on the surface, its Riemannian
metric and its spin structure.

Let us show that the results for closed surfaces carry over to the 
complete noncompact case without any essential changes.

\begin{theorem}\label{ghochnc}
Let $M$ be a complete surface of genus $g\ge 1$ with a Riemannian metric
of finite area.
Let $M$ be equipped with a spin structure which is nontrivial along the
ends and which admits a spin-cut.
Let $\de(M)$ be the  spin-cut-diameter of $M$.
Then
$$
\la_\ast \geq {2\sqrt{\pi}\over (2g+1)\,\sqrt{\area(M)}} -{1\over \de(M)}.
$$
If $g=1$, then for any $k\in \NN$
$$
\la_\ast\geq -{2\over k\,\de(M)} + 
\sqrt{{\pi\over k\,\area(T^2)}+{2\over k^2 \de(M)^2}}.
$$
\end{theorem}

\proof{}
Let $\ep >0$ and let $\ga_1, \ldots ,\ga_g$ be a spin-cut such that its 
spin-cut-diameter satisfies
$$
\de(\ga_1,\ldots,\ga_g) \ge \de(M) - \ep.
$$
Pick a smooth spinor $\ph$ on $M$ with compact support such that
$$
\frac{\|D\ph\|^2}{\|\ph\|^2} \le \la_\ast + \ep.
$$
Now we change the metric on $M$ outside the support of $\ph$ and away from
the $\ga_i$ such that it extends to $\barM$ and such that
$$
\area(\barM) \le \area(M) + \ep.
$$
Since the spin structure of $M$ is nontrivial along the ends it extends to one
on~$\barM$.
Theorem~\ref{ghoch} applied to $\barM$ now yields
\begin{eqnarray*}
\la_\ast + \ep& \ge & \frac{\|D\ph\|^2}{\|\ph\|^2} \\
&\ge& {2\sqrt{\pi}\over (2g+1)\,\sqrt{\area(\barM)}} -
{1\over \de(\ga_1,\ldots,\ga_g)} \\
&\ge&
{2\sqrt{\pi}\over (2g+1)\,\sqrt{\area(M)+\ep}} -
{1\over \de(M) - \ep}.
\end{eqnarray*}
Taking $\ep\to 0$ finishes the proof of the first assertion.
The second part for $g=1$ is shown similarly.
\qed

The assumption that the spin structure be nontrivial along the ends is crucial.
It has been shown by the second author \cite{baer:p98} that the $L^2$-spectrum
of the Dirac operator on a complete hyperbolic surface of finite area
whose spin structure is not nontrivial along the ends is given by
  $$
  \lss(D) = \RR.
  $$

%%%%%%%%%%%%%%%%%%%%%%%%%%%%%%%%%%%
\appendix

%%%%%%%%%%%%%%%%%%%%%%%%%%%%%%%%%%%
\section{Two lemmata about cylinders}
%%%%%%%%%%%%%%%%%%%%%%%%%%%%%%%%%%%

\begin{lemma}\label{appendixlemmaA}
Let $\ga:S^1\to S^2\setminus\{N,S\}$ be a simple closed curve in the 2-sphere
without North Pole $N$ and South Pole $S$. 
Then either $\ga$ is contractible in $S^2\setminus\{N,S\}$
or the homotopy class of $\ga$ generates 
$\pi_1(S^2\setminus\{N,S\})\cong \ZZ$.
\end{lemma}

\proof{}
According to the theorem of Jordan-Schoenfliess there is a 
diffeomorphism $\ph:S^2\to S^2$ mapping $\ga$ to the equator.
If $\phi(S)$ and $\phi(N)$ lie in the same hemisphere, then 
$\ga$ bounds a disk in $Z=S^2\setminus\{N,S\}$. 
In this case $\ga$ is contractible in $S^2\setminus\{N,S\}$.
Otherwise $[\ga]$ generates the fundamental group of $S^2\setminus\{N,S\}$.
\qed

\begin{lemma}\label{appendixlemmaB}
Let $Z:=\left\{(x,y,z)\,|\,x^2+y^2=1\right\}\subset \RR^3$ be the cylinder.
Let $f:Z\to \RR$ be smooth and assume that
$f(x,y,z)\to \infty$ for $z\to \infty$ and $f(x,y,z)\to -\infty$ for $z\to -\infty$ uniformly in $x,y$. 
This is equivalent to assuming that $f$ is proper and onto.
Then for any regular value $t\in\RR$ the set $f^{-1}(t)$ has a connected 
component which is a simple closed curve whose homotopy class generates
$\pi_1(Z)$.
\end{lemma}

\proof{}
Since $f$ is proper and $t$ is regular $N:=f^{-1}(t)$ is a closed 
$1$-dimensional manifold, \ie a finite union of simple closed curves.
Not every connected component of $N$ is contractible in $Z$, as otherwise
for large $K$ it would be possible to connect
$(1,0,-K)$ and $(1,0,K)$ by a curve in $Z\setminus N$. This is impossible 
by the mean value theorem.

Let $\ga$ by a parametrization of a noncontractible component of $N$. 
According to the previous lemma $[\ga]$ generates $\pi_1(Z)$.
\qed

%%%%%%%%%%%%%%%%%%%%%%%%%%%%%%%%%%%%%%%%%%%%%%%%%%%%%%%%%%%%%%%%%%%%%%%%%%

\providecommand{\bysame}{\leavevmode\hbox to3em{\hrulefill}\thinspace}

\vspace{1cm}

\parskip0ex

Fachbereich Mathematik

Universit\"at Hamburg

Bundesstra\ss{}e 55

20146 Hamburg

Germany

\vspace{0.5cm}

E-Mail:
{\tt ammann@math.uni-hamburg.de}

\hspace{1.2cm}
{\tt baer@math.uni-hamburg.de}

WWW:
{\tt http://www.math.uni-hamburg.de/home/ammann/}

\hspace{1.18cm}
{\tt http://www.math.uni-hamburg.de/home/baer/}


\begin{thebibliography}{10}
\bibitem{ammann:p00a}
  B.~Ammann, \emph{Spectral estimates on 2-tori}, {P}reprint {A}pril 2000, 
  Hamburger Beitr\"age zur Mathematik,
  no.~95.

\bibitem{ammann:p00b}
\bysame, \emph{A spin-conformal lower bound of the first positive {D}irac
  eigenvalue}, {P}reprint {O}ctober 2000, Hamburger Beitr\"age zur
  Mathematik, no.~96.

\bibitem{bangert:90}
V.~Bangert, \emph{Minimal geodesics}, Ergodic Theory Dynamical Systems
  \textbf{10} (1990), 263--286.

\bibitem{baer:92a}
C.~B{\"a}r, \emph{The {D}irac operator on homogeneous spaces and its
  spectrum on 3-dimensional lens spaces}, Arch. Math. \textbf{59} (1992),
  65--79.

\bibitem{baer:92b}
\bysame, \emph{Lower eigenvalue estimates for {D}irac operators}, Math. Ann.
  \textbf{293} (1992), 39--46.

\bibitem{baer:94}
\bysame, \emph{The {D}irac operator on space forms of positive curvature}, J.
  Math. Soc. Japan \textbf{48} (1994), 69--83.

\bibitem{baer:96}
\bysame, \emph{Metrics with harmonic spinors}, Geom. Funct. Anal. \textbf{6}
  (1996), 899--942.

\bibitem{baer:p98}
\bysame, \emph{The {D}irac operator on hyperbolic manifolds of finite
  volume}, J. Diff. Geom. \textbf{54} (2000), 439--488.

\bibitem{baer:98}
\bysame, \emph{Extrinsic bounds for eigenvalues of the {D}irac operator}, Ann.
  Global Anal. Geom. \textbf{16} (1998), 573--596.

\bibitem{bunke:91}
U.~Bunke, \emph{Upper bounds of small eigenvalues of the {D}irac operator and
  isometric immersions}, Ann. Glob. Anal. Geom. \textbf{9} (1991), 109--116.

\bibitem{burago:92}
D.~Yu. Burago, \emph{Periodic metrics}, In: Representation theory and dynamical
  systems, Amer. Math. Soc., Providence, RI, 1992, 205--210.

\bibitem{cahen.franc.gutt:89}
M.~Cahen, A.~Franc, and S.~Gutt, \emph{Spectrum of the {D}irac operator on
  complex projective space ${P}_{2q-1}({\CC})$}, Lett. Math. Phys. \textbf{18}
  (1989), 165--176.

\bibitem{cahen.franc.gutt:94}
\bysame, \emph{Erratum to 'spectrum of the {D}irac operator on complex
  projective space ${P}_{2q-1}({\CC})$'}, Lett. Math. Phys. \textbf{32} (1994),
  365--368.

\bibitem{camporesi.higuchi:96}
R.~Camporesi and A.~Higuchi, \emph{On the eigenfunctions of the
  {D}irac operator on spheres and real hyperbolic spaces}, J. Geom. Phys.
  \textbf{20} (1996), 1--18.

\bibitem{federer:96}
H.~Federer, \emph{Geometric measure theory}, Springer-Verlag New York, 1969.

\bibitem{fegan:87}
H.~Fegan, \emph{The spectrum of the {D}irac operator on a simply connected
  compact {L}ie group}, Simon Stevin \textbf{61} (1987), 97--108.

\bibitem{friedrich:80}
T.~Friedrich, \emph{Der erste {E}igenwert des {D}irac-{O}perators einer
  kompakten {R}iemannschen {M}annigfaltigkeit nicht-negativer {K}r{\"u}mmung},
  Math. Nach. \textbf{97} (1980), 117--146.

\bibitem{friedrich:84}
\bysame, \emph{Zur {A}bh{\"a}ngigkeit des {D}irac-{O}perators von der
  {S}pin-{S}truk\-tur}, Colloq. Math. \textbf{48} (1984), 57--62.

\bibitem{friedrich:buch}
\bysame, \emph{Dirac operators in Riemannian geometry},
  Graduate Studies in Mathematics. 25. Providence, 
  RI: American Mathematical Society (2000)

\bibitem{gromov:81b}
M.~Gromov, \emph{Groups of polynomial growth and expanding maps}, Inst.
  Hautes \'Etudes Sci. Publ. Math. \textbf{53} (1981), 53--73.

\bibitem{gromov:99}
\bysame, \emph{Metric structures for {R}iemannian and non-{R}iemannian
  spaces}, Birkh\"auser Boston, 1999.

\bibitem{hijazi:86}
O.~Hijazi, \emph{A conformal lower bound for the smallest eigenvalue of
  the {D}irac operator and {K}illing spinors}, Comm. Math. Phys. \textbf{104}
  (1986), 151--162.

\bibitem{hijazi:91}
\bysame, \emph{Premi\`ere valeur propre de l'op\'erateur de {D}irac et nombre
  de {Y}amabe}, C. R. Acad. Sci. Paris \textbf{t.\ 313, S\'erie I} (1991),
  865--868.

\bibitem{hitchin:74}
N.~Hitchin, \emph{Harmonic spinors}, Adv. Math. \textbf{14} (1974), 1--55.

\bibitem{kirchberg:86}
K.-D. Kirchberg, \emph{An estimation for the first eigenvalue of the {D}irac
  operator on closed {K}\"ahler manifolds of positive scalar curvature}, Ann.
  Glob. Anal. Geom. \textbf{4} (1986), 291--325.

\bibitem{kirchberg:88}
\bysame, \emph{Compact six-dimensional {K}{\"a}hler spin manifolds of positive
  scalar curvature with the smallest possible first eigenvalue of the {D}irac
  operator}, Math.\ Ann. \textbf{282} (1988), 157--176.

\bibitem{kramer.semmelmann.weingart:98}
W.~Kramer, U.~Semmelmann, and G.~Weingart, \emph{The first eigenvalue of the
  {D}irac operator on quaternionic {K}\"ahler manifolds}, Comm. Math. Phys.
  \textbf{199} (1998), 327--349.

\bibitem{kramer.semmelmann.weingart:99}
\bysame, \emph{Eigenvalue estimates for the {D}irac operator on quaternionic
  {K}\"ahler manifolds}, Math. Z. \textbf{230} (1999), 727--751.

\bibitem{kusner.schmitt:p96}
R.~Kusner and N.~Schmitt, \emph{The spinor representation of surfaces in
  space}, Preprint, 1996, http://www.arxiv.org/abs/dg-ga/9610005.

\bibitem{lawson.michelsohn:89}
H.-B.~Lawson and M.-L. Michelsohn, \emph{Spin geometry}, Princeton University
  Press, Princeton, 1989.

\bibitem{lott:86}
J.~Lott, \emph{Eigenvalue bounds for the {D}irac operator}, Pacific J. of
  Math. \textbf{125} (1986), 117--126.

\bibitem{massart:97a}
D.~Massart, \emph{{Stable norms of surfaces: local structure of the unit ball
  at rational directions.}}, Geom. Funct. Anal. \textbf{7} (1997),
  996--1010.

\bibitem{milhorat:p}
J.-L. Milhorat, \emph{Spectrum of the Dirac operator on ${G}r_2({\CC}^{m+2})$},
  J. Math. Phys. \textbf{39} (1998), 594-609.

\bibitem{milhorat:92}
\bysame, \emph{Spectre de l'op\'erateur de {D}irac sur les espaces projectifs
  quaternioniens}, C. R. Acad. Sci. Paris \textbf{1} (1992), 69--72.

\bibitem{pansu:83}
P.~Pansu, \emph{Croissance des boules et des g\'eod\'esiques ferm\'ees dans
  les nilvari\'et\'es}, Ergodic Theory Dyn. Syst. \textbf{3} (1983),
  415--445.

\bibitem{pfaeffle:p99}
F.~Pf{\"a}ffle, \emph{The {D}irac spectrum of {B}ieberbach manifolds},
  J. Geom. Phys. \textbf{35} (2000), 367--385.

\bibitem{pinkall:85a}
U.~Pinkall, \emph{Regular homotopy classes of immersed surfaces}, Topology
  \textbf{24} (1985), 421--434.

\bibitem{roe:88}
J.~Roe, \emph{Elliptic operators, topology and asymptotic methods}, Pitman
  Research Notes in Mathematics Series, no. 179, Longman, 1988.

\bibitem{seegerdiplom}
L.~Seeger, \emph{{D}er {D}irac-{O}perator auf kompakten symmetrischen
  {R}{\"a}umen}, Diplomarbeit, Universit{\"a}t Bonn, 1997.

\bibitem{seeger:99}
\bysame, \emph{The spectrum of the {D}irac operator on ${G}\sb 2/{\rm
  {S}{O}}(4)$}, Ann. Global Anal. Geom. \textbf{17} (1999), 385--396.

\bibitem{seegerdiss}
\bysame, \emph{Metriken mit harmonischen {S}pinoren auf geradedimensionalen 
  {S}ph{\"a}ren}, Dissertation, {U}niversit\"at {H}amburg, 2000.

\bibitem{seifarth.semmelmann:93}
S.~Seifarth and U.~Semmelmann, \emph{The spectrum of the {D}irac operator on
  the odd dimensional complex projective space ${P}^{2m-1}({\CC})$},
  SFB288-Preprint~95, 1993.

\bibitem{strese:80b}
H.~Strese, \emph{Spektren symmetrischer {R}{\"a}ume}, Math. Nachr. \textbf{98}
  (1980), 75--82.

\bibitem{strese:80a}
\bysame, \emph{\"uber den {D}irac-{O}perator auf {G}ra\ss
  mann-{M}annigfaltigkeiten}, Math. Nachr. \textbf{98} (1980), 53--59.

\bibitem{sulanke:79}
  S.~Sulanke, \emph{{D}ie {B}erechnung des {S}pektrums des {Q}uadrates des
  {D}irac-{O}perators auf der {Sp}h{\"a}re}, Dissertation,
  Humboldt-Universit{\"a}t, Berlin, 1979.

\bibitem{topping:00} 
  P.~Topping, \emph{Towards the {W}illmore conjecture}, 
  Calc. Var. Partial Differential Equations \textbf{11} (2000), 361--393.

\bibitem{trautman:93}
A.~Trautman, \emph{Spin structures on hypersurfaces and the spectrum of the
  {D}irac operator on spheres}, In: Oziewicz, Zbigniew (ed.) et al., 
  Spinors, Twistors, {C}lifford Algebras and
  Quantum Deformations, Kluwer Academic Publishers, 1993, 25-29.

\bibitem{trautman:95}
\bysame, \emph{The {D}irac operator on hypersurfaces}, Acta Phys.~Polon.~B
  \textbf{26} (1995), 1283--1310.

\end{thebibliography}
\end{document}